\documentclass[a4paper]{article}
\usepackage[utf8]{inputenc} % Character encoding
\usepackage{biblatex} % Imports biblatex package
\usepackage{amsmath, amssymb, amsthm, enumerate}
\allowdisplaybreaks
\usepackage{xcolor} % Colours
\usepackage{todonotes}
\usepackage[a4paper, margin=3cm]{geometry}
\usepackage{xspace}
\usepackage{tikz}
\usepackage{mathtools}
\usepackage{mathtools}
\usepackage{esint}

\usepackage{enumitem}
\usepackage{stackengine}

\DeclareMathOperator{\supp}{supp}

\newcommand\preceqdot{\mathrel{\ooalign{$\prec$\cr
  \hidewidth\raise0ex\hbox{$\cdot\mkern-0.5mu$}\cr}}}

\newcommand{\R}{{\bf R}}
\newcommand{\C}{{\bf C}}
\newcommand{\Z}{{\bf Z}}
\newcommand{\N}{{\bf N}}

\DeclareMathOperator{\dist}{dist}
\newtheorem{thm}{Theorem}[section]
\newtheorem{lem}[thm]{Lemma}

\newtheorem{proposition}[thm]{Proposition}
{\theoremstyle{definition}
\newtheorem{definition}[thm]{Definition}

\newtheorem{remark}[thm]{Remark}
}
\title{Some More Sparse Bounds for Rough and Smooth Pseudodifferential Operators}

\author{Solange Mukeshimana \and David Rule}

\date{25th March 2026}

\newcommand{\Addresses}{{% additional braces for segregating \footnotesize
  \bigskip
  \footnotesize

  S.~Mukeshimana, \textsc{College of Science and Technology, University of Rwanda, P.O.~Box: 3900, Kigali, Rwanda}\par\nopagebreak
  \textit{E-mail address}: \texttt{sosmukish@gmail.com}

  \medskip

  D.~Rule (Corresponding author), \textsc{Department of Mathematics, Linköping University, SE-581 83 Linköping, Sweden}\par\nopagebreak
  \textit{E-mail address}: \texttt{david.rule@liu.se}

}}

\begin{document}

\maketitle

\begin{abstract}
Beltran \& Cladek~\cite{BC} use $L^r$ to $L^s$ bounds to prove sparse form bounds for pseudodifferential operators with H\"ormander symbols in $S^m_{\rho,\delta}$ up to, but not including, the sharp end-point in decay $m$. We further develop their technique, obtaining pointwise sparse bounds for rough pseudodifferential operators that are merely measurable in their spatial variables and an alternative proof of their results which avoids proving geometrically decaying sparse bounds. We also provide sufficient conditions for sparse form bounds to hold and use these to reprove know sparse bounds for pseudodifferential operators with symbols in $S^0_{1,\delta}$ for $\delta < 1$.
\end{abstract}

%%%%%%%%%%%%%%%%%%%%%%
%%%%%%%%%%%%%%%%%%%%%%
%%%%%%%%%%%%%%%%%%%%%%

\section{Introduction and Main Results}\label{introduction}

In this paper, we investigate the boundedness properties of various pseudodifferential operators of H\"ormander type. Given a function $a(x,y,\xi)$ (often referred to as an \emph{amplitude}) we define the \emph{pseudodifferential operator} $T_a$ by
\begin{equation*}
    T_a(f)(x) = \int_{\R^n} \int_{\R^n} a(x,y,\xi) f(y) e^{i \xi\cdot(x-y)} dy d\xi.
\end{equation*}
In the case that $a = a(x,\xi)$ does not depend on $y$ (and as such is often referred to as a \emph{symbol}) the operator can be written as
\begin{equation*}
    T_a(f)(x) = \int_{\R^n} a(x,\xi) \widehat{f}(\xi) e^{i \xi\cdot x} d\xi.
\end{equation*}
First, we recall the H\"ormander classes of symbols.
\begin{definition}\label{def:symbol}
     A function $a \colon \R^n \times \R^n \to \C$ belongs to the set $S^m_{\rho,\delta}$ ($m,\rho,\delta\in\R$) if it is infinitely differentiable and satisfies
    \begin{equation*}
        |\partial^\alpha_\xi\partial^\beta_x a(x,\xi)| \lesssim (1+|\xi|)^{m-\rho|\alpha|+\delta|\beta|}
    \end{equation*}
    for each pair of multi-indices $\alpha$ and $\beta$.
\end{definition}
We also recall the definition of rough symbols and amplitudes (see, for example, \cite{MRS1}).

\begin{definition}\label{def:rough_amplitude}
     A function $a \colon \R^n \times \R^n \times \R^n \to \C$ is said to belong to the set $L^\infty A^m_\rho$ ($m,\rho\in\R$) if
    \begin{equation*}
        |\partial^\alpha_\xi a(x,y,\xi)| \lesssim (1+|\xi|)^{m-\rho|\alpha|}
    \end{equation*}
    for each multi-index $\alpha$ and is measurable in the $x$ and $y$ variables.
\end{definition}
\begin{definition}\label{def:rough_symbol}
    A function $a \colon \R^n \times \R^n \to \C$ is said to belong to the set $L^\infty S^m_\rho$ if
    \begin{equation*}
        |\partial^\alpha_\xi a(x,\xi)| \lesssim (1+|\xi|)^{m-\rho|\alpha|}
    \end{equation*}
    for each multi-index $\alpha$, and is measurable in the $x$ variable.
\end{definition}

Beltran \& Cladek~\cite{BC} proved that for $a\in S^m_{\rho,\delta}$, with $0 < \delta \leq \rho < 1$ and
\begin{equation}\label{limitingvalues}
    m < -n(1-\rho)\left(\frac{1}{r} - \frac{1}{s}\right), \quad \mbox{for $1\leq r\leq 2 \leq s \leq \infty$,}
\end{equation}
there exists a constant $C=C(m,\rho,\delta,r,s)$ such that for each pair of bounded functions $f$ and $g$ with compact support, there exists a sparse collection $\mathcal{S}$ such that
\begin{equation} \label{ineq:sparseform}
    \left|\langle T_a(f),g\rangle\right| \leq C \sum_{Q\in\mathcal{S}} \langle f \rangle_{r,Q}\langle g \rangle_{s',Q}|Q|.
\end{equation}
A collection $\mathcal{S}$ of sets is said to be \emph{$\eta$-sparse} for some $\eta \in (0,1]$, or for brevity simply \emph{sparse}, if, for each $Q \in \mathcal{S}$, there exists a set $E(Q)$ such that $E(Q) \subseteq Q$ and $|E(Q)| \geq \eta|Q|$, and that the collection $\{E(Q) \colon Q\in\mathcal{S}\}$ is pairwise disjoint. For $r < \infty$, the notation
    \begin{equation*}
        \langle f \rangle_{r,Q} := \left(\frac{1}{|Q|} \int_Q |f|^r \right)^\frac{1}{r}
    \end{equation*}
is the $L^r$-average over a set $Q$ and $\langle f \rangle_{\infty,Q} := \sup_Q |f|$. Estimate \eqref{ineq:sparseform} is called a \emph{sparse form bound} and there are many well-known consequences of such bounds (see, for example, \cite{BC}, \cite{BernicotFreyPetermichl}, \cite{Conde-AlonsoCuliucPlinioOu} and \cite{LiPerezRivera-RiosRoncal}). Here and throughout the paper, for an exponent $s$, $s'$ will denote the dual exponent of $s$, so $1/s + 1/s' = 1$. The sparse form bound \eqref{ineq:sparseform} immediately extends to larger values of $r$ and $s'$ outside the range $0\leq r\leq 2 \leq s < \infty$, since H\"older's inequality gives that averages $\langle f \rangle_{r,Q}$ are monotone increasing in the exponent $r$.

We can encapsulate the method in \cite{BC} in the following proposition, which says that a decomposition of the operator into localised pieces, plus a scale-invariant $L^{r}$ to $L^{s}$ bound on those localised pieces leads to sparse form bounds. Here and throughout the paper, $\frac{1}{3}Q$ denotes a cube $Q$ scaled by a factor $\frac{1}{3}$ about its centre.
\begin{proposition}\label{intro:sparseform}
Let $Q_k$ denote a cube of radius $2^{-k}$ with $k\in\Z$ and assume that the operator $T$ is a countable sum $T = \sum_{j\in I} T^{j}$ of sublinear operators $T_j$ indexed by $j \in I$. Furthermore, assume that, for a given function $\iota \colon I \to \Z$, exponents $1\leq r,s \leq\infty$ and constants $B_j$ ($j\in I$) such that $\sum_{j\in I} B_{j} < \infty$, we have the estimate
\begin{equation}\label{intro:Lr_Ls_averagebounds}
    \left\langle T^{j}\left(f\chi_{\frac{1}{3}Q_{\iota(j)}}\right)\right\rangle_{s,Q_{\iota(j)}} \leq B_{j} \left\langle f \right\rangle_{r,Q_{\iota(j)}}
\end{equation}
for each $j \in I$ and $T^{j}\left(f\chi_{\frac{1}{3}Q_{\iota(j)}}\right)$ is supported in $Q_{\iota(j)}$. Then, there exists a constant $C$, such that for each pair of  bounded and compactly supported functions $f$ and $g$, there exists a sparse collection $\mathcal{S}$ such that
\begin{equation*}
    \left\langle T(f), g\right\rangle \leq C \sum_{Q \in \mathcal{S}} \langle f \rangle_{r,Q}\langle g \rangle_{s',Q}|Q|.
\end{equation*}
\end{proposition}
We note that the averages in \eqref{intro:Lr_Ls_averagebounds} are averages on the cube $Q_{\iota(j)}$ and so include scaling factors of $|Q_{\iota(j)}|^{-1/s}$ and $|Q_{\iota(j)}|^{-1/r}$, respectively, on each side of the inequality. We deliberately choose to write the inequality in terms of averages instead of simple $L^p$-norms as it is the scale-invariant average version that must be summable over the pieces of the operator.

We provide an alternative proof of Proposition~\ref{intro:sparseform} as this will allow us to extend the methods in \cite{BC} to obtain \emph{pointwise sparse bounds} (that is, the estimate \eqref{def:pwsparse} below) in the case where \eqref{intro:Lr_Ls_averagebounds} holds with $s=\infty$.
\begin{proposition}\label{intro:pwsparse}
Let $Q_k$ denote a cube of radius $2^{-k}$ with $k\in\Z$ and assume that the operator $T$ is a countable sum $T = \sum_{j\in I} T^{j}$ of sublinear operators $T_j$ indexed by $j \in I$. Furthermore, assume that, for a given function $m \colon I \to \Z$, exponent $1\leq r \leq \infty$ and constants $B_j$ ($j\in I$) such that $\sum_{j\in I} B_{j} < \infty$, we have the estimate
\begin{equation}\label{intro:pwparseaveragebound}
    \left|T^{j}\left(f\chi_{\frac{1}{3}Q_{\iota(j)}}\right)(x)\right| \leq B_{j} \left\langle f \right\rangle_{r,Q_{\iota(j)}}\chi_{Q_\iota(j)}(x).
\end{equation}
Then there exists a constant $C$, such that for each bounded and compactly supported function $f$ there exists a sparse collection $\mathcal{S}$ such that
\begin{equation}\label{def:pwsparse}
    \left|T(f)\right| \leq C \sum_{Q \in \mathcal{S}} \langle f \rangle_{r,Q}\chi_Q(x).
\end{equation}
\end{proposition}

We apply this proposition to extend the results in \cite{BC} in a few ways. When $s=\infty$ in \eqref{limitingvalues} we obtain pointwise sparse bounds (with exponent $r$) rather than sparse form bounds. Moreover, we obtain them for a wider class of operators, namely those arising from symbols in $L^\infty S^m_\rho$ when $r\geq1$ and from amplitudes in $L^\infty A^m_\rho$ when $r=1$. As such, they are also generalisations of Proposition~3.7.3 in \cite{BeltranPortales}. We formulate these extensions as the following theorems.

\begin{thm}\label{thm:rough-symbols}
    Suppose $a\in L^\infty S^m_{\rho}$, with $\rho \in [0,1]$,
\begin{equation*}
    m < -\frac{n}{r}(1-\rho)
\end{equation*}
and $r \in [1,2]$. Then there exists a constant $C$, which depends on $m$, $\rho$, $r$ and a finite number of the implicit constants in Definition~\ref{def:symbol}, such that for each bounded function $f$ with compact support, there exists a sparse collection $\mathcal{S}$ such that
\begin{equation*}
    \left|T_a(f)(x)\right| \leq C \sum_{Q \in \mathcal{S}} \langle f \rangle_{r,Q}\chi_Q(x).
\end{equation*}
\end{thm}
\begin{thm}\label{thm:rough-amplitudes}
   Suppose $a\in L^\infty A^{m}_{\rho}$ with $\rho \in [0,1]$ and $m < -n(1-\rho)$. Then there exists a constant $C$, which depends on $m$, $\rho$ and a finite number of the implicit constants in Definition~\ref{def:rough_amplitude}, such that for each bounded function $f$ with compact support, there exists a sparse collection $\mathcal{S}$ such that
\begin{equation*}
    \left|T_a(f)(x)\right| \leq C \sum_{Q \in \mathcal{S}} \langle f \rangle_{1,Q}\chi_Q(x).
\end{equation*}
\end{thm}

We note that the methods in \cite{BC} also produce sparse form bounds in some special cases of the weaker hypotheses appearing in Theorems~\ref{thm:rough-symbols} and \ref{thm:rough-amplitudes}. Specifically, they do not use $x$-smoothness to prove \eqref{intro:Lr_Ls_averagebounds} for $r=1$ or $r=\infty$ and $s=\infty$. They only use it to allow adjoints of operators to be dealt with via the same techniques and for proving \eqref{intro:Lr_Ls_averagebounds} in the case $r=s=2$. To obtain the full versions of Theorems~\ref{thm:rough-symbols} and \ref{thm:rough-amplitudes}, we replace the use of adjoints (with their implicit use of $x$-smoothness) with a more abstract $TT^*$ argument. Beltran \& Cladek also prove some pointwise sparse bounds for operators with smooth symbols, which they formulate as Theorem~1.3 in \cite{BC}, although they use a different approach to the methods of Proposition~\ref{intro:pwsparse} below. More specifically, they use a sharp maximal function estimate (see Definition~\ref{def:maximal} below) from \cite{Michalowski_Rule_Staubach_2012}~--- which requires smoothness in $x$~--- and Lerner's pointwise domination principle~\cite{Lerner2016}. In the case $\rho < 1$, Theorem~\ref{thm:rough-symbols} is more general than Theorem~1.3 in \cite{BC} as it admits a wider range of $m$ and a larger class of symbols.
\begin{remark}\label{commentc}
    The methods we present here provide an alternative route to proving weak-type estimates for pseudodifferential operators. Conde-Alonso, Culiuc, Di Plinio \& Ou~\cite{Conde-AlonsoCuliucPlinioOu} prove that if a sparse form bound \eqref{ineq:sparseform} holds for $r=1$ and some $s' \in [1,\infty)$, then there exists a constant $C>0$ such that the weak-type estimate
    \begin{equation}\label{ineq:weak-type}
        \left|\{x \colon |T_a(f)(x)|>\lambda\}\right| \leq \frac{C}{\lambda} \int |f| 
    \end{equation}
    holds. Since (in contrast to \cite{BeltranPortales}) we do not use the weak-type boundedness of pseudodifferential operators to obtain the sparse bounds of Theorems~\ref{thm:rough-symbols} and \ref{thm:rough-amplitudes}, it follows from Theorem~\ref{thm:rough-amplitudes} above and Theorem~E in \cite{Conde-AlonsoCuliucPlinioOu} that \eqref{ineq:weak-type} holds when $a \in L^\infty A^m_\rho$ for $m < -n(1-\rho)$. Moreover, the argument can be pushed a little further in the case of smooth symbols, although it just fails to reach the end-point value of $m=-n(1-\rho)/2$ of Theorem~3.2 in \cite{AlvarezHounie90} for $\delta<\rho$. To see this, recall that the collection of operators $T_a$ for $a\in S^m_{\rho,\delta}$ and $\delta < \rho$ is closed under the operation of taking adjoints. Therefore, for a given $a\in S^m_{\rho,\delta}$ with $\delta<\rho$, there exists a $b\in S^m_{\rho,\delta}$, such that $T_b = (T_a)^*$. We can then calculate that
    \begin{equation*}
        \langle T_a(f), g \rangle = \langle f, (T_a)^*(g) \rangle = \langle f, T_b(g)\rangle
        \lesssim \sum_{Q\in\mathcal{S}} \langle f \rangle_{1,Q}\langle g \rangle_{2,Q}|Q|,
    \end{equation*}
    where the last estimate follows by applying Theorem~\ref{thm:rough-symbols} with $a=b$ and $r=2$. Thus, the hypotheses of Theorem~E in \cite{Conde-AlonsoCuliucPlinioOu} are fulfilled and \eqref{ineq:weak-type} follows for $a\in S^m_{\rho,\delta}$ when $m < -n(1-\rho)/2$ and $\delta < \rho$.
\end{remark}

A natural next step to take is to ask: can we not further develop these techniques for sparse form bounds to obtain results at the end-point $m = n(\rho-1)\left(\frac{1}{r} - \frac{1}{s}\right)$ for smooth symbols? To discuss this further we introduce some relevant maximal functions.
\begin{definition} \label{def:maximal}
    We denote the usual uncentred Hardy-Littlewood maximal operator on balls by $M$ and, for $p\geq1$, the $L^p$-maximal operator by
    \begin{equation}\label{def:lphlmaximal}
        M_p(f)(x)  = \sup_{B\ni x} \left(\frac{1}{|B|} \int_B |f|^p\right)^{\frac{1}{p}},
    \end{equation}
    so $M = M_1$. We write $M^c_p$ to mean the uncentred maximal function with respect to cubes, more precisely, \eqref{def:lphlmaximal} with the supremum over balls $B$ replaced by that over cubes $Q$ with sides parallel to the axes. We define the sharp maximal function as
    \begin{equation*}
        M^\sharp(f)(x)  = \sup_{B\ni x} \frac{1}{|B|} \int_B |f(y) - f_{B}|dy,
    \end{equation*}
    where $f_B = \frac{1}{|B|}\int_B f$.
\end{definition}
Before the theory of sparse bounds was developed, a typical route to proving, for example, weighted norm inequalities (albeit without the same good control of how the constants depend on the weight, as is obtained with sparse bounds) would be to prove the pointwise inequality
\begin{equation}\label{sharpfunctionestimate}
M^\sharp \circ T_a(f)(x) \lesssim M_{r}(f)(x).
\end{equation}
These methods were pioneered by Chanillo \& Torchinsky~\cite{ChanilloTorchinsky1986} (where $r=2$, $m\leq -n(1-\rho)/2$ and $\delta < \rho$) and Alvarez \& Hounie~\cite{AlvarezHounie90} (where $r > 1$ and $m\leq -n(1-\rho)$). And even more recently, Wang~\cite{wang2022sharpfunctionweightedlp} has improved on these results (where $r=2$, $m\leq -n(1-\rho)/2$~--- see also Chen \& Wang~\cite{ChenWang2023}).

In an attempt to answer the question of what can be said at the end-point $m = n(\rho-1)\left(\frac{1}{r} - \frac{1}{s}\right)$, we try to push the ideas of Propositions~\ref{intro:sparseform} and \ref{intro:pwsparse} further. The result of this thought process will be formulated later in the article in Proposition~\ref{intro:end-pointsparsebdd2}. As we cannot expect to sum pieces of an operator at the end-point, we instead use an induction argument over cubes (analogous to \cite{LMSvideos}) and the hypothesis \eqref{intro:Lr_Ls_averagebounds} is replaced by
\begin{equation*}
        \displaystyle\left\langle T(f\chi_{Q\setminus Q'}) - (T(f\chi_{Q\setminus Q'}))_{\frac{1}{3}Q'} \right\rangle_{s,\frac{1}{3}Q'} \leq A_1 \left\langle f \right\rangle_{r,Q}.
\end{equation*}
where $Q' \subseteq Q$ are cubes such that $Q$ is on the preceding step to $Q'$ in the induction argument and $(T(f\chi_{Q\setminus Q'}))_{\frac{1}{3}Q'} = |\frac{1}{3}Q'|^{-1}\int_{\frac{1}{3}Q'} T(f\chi_{Q\setminus Q'})$. We apply this proposition to pseudodifferential operators in Theorem~\ref{thm:endpoint-pseudos} below, although frustratingly, we only succeed for $\rho=1$. As such, we manage to prove a sparse form version of Theorem~1.3 from \cite{BC} in the case $\rho=1$. As they note, this result was already known with the expected $L^1$-averages via the theory of Calder\'on-Zygmund operators. We hope that future work will broaden the applicability of Proposition~\ref{intro:end-pointsparsebdd2}. As a by-product of this work, we succeed in proving \eqref{sharpfunctionestimate} under weaker hypotheses than previously known when $r < 2$, which we formulate as Theorem~\ref{sharpestimatethm} below.
\begin{thm} \label{sharpestimatethm}
Suppose $\lambda := \max\{0,(\delta-\rho)/2\} < 1/n$ and let $a\in S^{m}_{\rho, \delta}$, with $0 <\rho\leq 1$, $0\leq \delta < 1$ and
\begin{equation*}
    m = -\frac{n}{p}(1-\rho) - n\lambda
\end{equation*}
for $p \in (1,2] \cap [2\rho,2]$. Then, for bounded functions $f$ with compact support, one has that
\begin{equation}\label{sharppwbound1}
M^\sharp \circ T_a(f)(x_0) \lesssim M_{p}(f)(x_0).
\end{equation}
\end{thm}

In Section~\ref{sec:2} we prove the abstract results: Propositions~\ref{intro:sparseform}, \ref{intro:pwsparse} and \ref{intro:end-pointsparsebdd2}. In Section~\ref{sec:3} we prove the results for rough pseudodifferential operators (Theorems~\ref{thm:rough-symbols} and \ref{thm:rough-amplitudes}) and in Section~\ref{sec:4} we prove the end-point result for smooth symbols (Theorem~\ref{thm:endpoint-pseudos}) and the pointwise maximal function estimate (Theorem~\ref{sharpestimatethm}).

%%%%%%%%%%%%%%%%%%%%%%
%%%%%%%%%%%%%%%%%%%%%%
%%%%%%%%%%%%%%%%%%%%%%

\section{Sparse Bounds from Scale-Invariant $L^r$ to $L^s$ Bounds}
\label{sec:2}

The aim of this section is to prove Propositions~\ref{intro:sparseform}, \ref{intro:pwsparse} and \ref{intro:end-pointsparsebdd2}, which all require some form of localisation and scale-invariant $L^r$ to $L^s$ bounds. The main task is to construct appropriate sparse collections of dyadic cubes from these assumptions.

A \emph{(shifted) dyadic cube} with side length $2^{-k}$ is a Cartesian product
\begin{equation}\label{def:dyadic_cube}
Q^{m,\omega}_k := \bigtimes_{i=1}^n\left[2^{-k}\left(m_i+\frac{\omega_i}{3}\right),2^{-k}\left(m_i+1+\frac{\omega_i}{3}\right)\right)
\end{equation}
of half-open intervals, where $k \in \Z$, $m=(m_1,\dots,m_n)$, $\omega = (\omega_1,\dots,\omega_n)$, and $m_i \in\Z$ and $\omega_i \in \{-1,0,1,\}$ for $i=1,\dots,n$. We denote by $\mathcal{D}^\omega_k$ the collection of all such $Q^{m,\omega}_k$ with a given side length $2^{-k}$ and translation $\omega$ and $\mathcal{D}^\omega = \cup_k \mathcal{D}^\omega_k$. As is usual practice, for a fixed $\omega$, if one dyadic cube is contained in the other, say $Q' \subseteq Q$ we say that $Q'$ is a \emph{descendant} of $Q$ and that $Q$ is an \emph{ancestor} of $Q'$. When, in addition to this inclusion, the side length of $Q$ is double that of $Q'$, we say that $Q'$ is a \emph{child} of $Q$ and that $Q$ is the \emph{parent} of $Q'$.

We give a proof of Proposition~\ref{intro:sparseform}. We remind the reader that the proposition is essentially implicitly contained in \cite{BC}. The proof we present here avoids the intermediate step of proving geometrically decaying sparse form bounds.

\begin{proof}[Proof of Proposition~\ref{intro:sparseform}]\label{page:constructionofcubes}
Since $f$ and $g$ are bounded and compactly supported, we can select a sufficiently negative index $k$ of cube side length $2^{-k}$ so that the support of $f$ and $g$ is contained in the $2^n$ cubes $Q^{m,\omega}_k$ with $m_i = 0$ or $-1$ for all $i=1,\dots,n$ and $\omega \in \{-1,0,1,\}^n$. Denote by $\mathcal{S}^\omega_0$ the collection of all such cubes $Q^{m,\omega}_k$.

For fixed $\omega$, we now construct $\mathcal{S}^\omega_p$. We do this recursively, first for positive integers $p$, and then for negative $p$. We begin with positive $p$. We do this by selecting a collection $\mathcal{S}(Q)$ of sub-cubes of $Q$ for each $Q \in \mathcal{S}^\omega_{p-1}$ and then defining
\begin{equation*}
    \mathcal{S}^{\omega}_{p} = \bigcup_{Q \in \mathcal{S}^{\omega}_{p-1}} \mathcal{S}(Q).
\end{equation*}
Consider the collections
\begin{equation*}
    \tilde{\mathcal{S}}_f(Q) := \left\{Q' \in \mathcal{D}^\omega \colon \mbox{$Q' \subseteq Q$ and $\langle f \rangle_{r,Q'} >  4^{1/r}$}\langle f \rangle_{r,Q}\right\}
\end{equation*}
and
\begin{equation*}
    \tilde{\mathcal{S}}_g(Q) := \left\{Q' \in \mathcal{D}^\omega \colon \mbox{$Q' \subseteq Q$ and $\langle g \rangle_{s',Q'} >  4^{1/s'}\langle g \rangle_{s',Q}$}\right\}
\end{equation*}
of dyadic sub-cubes of $Q$ and, with $\tilde{\mathcal{S}}(Q) := \tilde{\mathcal{S}}_f(Q)\cup\tilde{\mathcal{S}}_g(Q)$, define
\begin{equation*}
    \mathcal{S}(Q) = \{Q' \in \tilde{\mathcal{S}}(Q) \colon \mbox{$Q' \not\subseteq Q''$ for any $Q'' \in \tilde{\mathcal{S}}(Q)$}\}
\end{equation*}
to be the collection of maximal cubes of $\tilde{\mathcal{S}}(Q)$ with respect to inclusion.

We define
\begin{equation*}
    E^c(Q) = \bigcup_{Q' \in \mathcal{S}(Q)} Q',
\end{equation*}
and $E(Q) = Q\setminus E^c(Q)$ for each dyadic cube $Q \in \mathcal{S}^\omega_p$. By construction, either
\begin{equation*}
    |Q'| \leq \left(\frac{1}{4^{1/r}\langle f \rangle_{r,Q}}\right)^r \int_{Q'} |f|^r \quad \mbox{or} \quad |Q'| \leq \left(\frac{1}{4^{1/s'}\langle g \rangle_{s',Q}}\right)^{s'} \int_{Q'} |g|^{s'}
\end{equation*}
for each $Q'\in\mathcal{S}(Q)$, and $\mathcal{S}(Q)$ is a disjoint union of cubes contained in $Q$. Therefore we can estimate
\begin{align*}
    |E^c(Q)| &= \sum_{Q' \in \mathcal{S}(Q)} |Q'| \\
    &\leq \sum_{Q' \in \mathcal{S}(Q)\cap\tilde{\mathcal{S}}_f(Q)} \left(\frac{1}{4^{1/r}\langle f \rangle_{r,Q}}\right)^r \int_{Q'} |f|^r + \sum_{Q' \in \mathcal{S}(Q)\cap\tilde{\mathcal{S}}_g(Q)} \left(\frac{1}{4^{1/s'}\langle g \rangle_{{s'},Q}}\right)^{s'} \int_{Q'} |g|^{s'} \\
    &\leq \left(\frac{1}{4^{1/r}\langle f \rangle_{r,Q}}\right)^r \int_{Q} |f|^r + \left(\frac{1}{4^{1/s'}\langle g \rangle_{{s'},Q}}\right)^{s'} \int_{Q} |g|^{s'} \\
    &\leq \left(\frac{\langle f\rangle_{r,Q}}{4^{1/r}\langle f \rangle_{r,Q}}\right)^r|Q| + \left(\frac{\langle g\rangle_{s',Q}}{4^{1/s'}\langle g \rangle_{s',Q}}\right)^{s'}|Q|
    \leq \frac{1}{2} |Q|.
\end{align*}
Thus, $|E(Q)| \geq (1-\frac{1}{2}) |Q| \geq \frac{1}{2}|Q|$ for each $Q \in \mathcal{S}^{\omega}_p$.

Now considering negative $p$, we define $\mathcal{S}^\omega_p$ to be the collection of (at most $2^n$) parents of cubes in $\mathcal{S}^\omega_{p+1}$ (with side lengths twice that of cubes in $\mathcal{S}^\omega_{p+1}$). Thus, for $Q \in \mathcal{S}^\omega_{p}$, there is a unique child $Q'$ of $Q$ contained in $\mathcal{S}^\omega_{p+1}$. With this choice of $Q$ and $Q'$, define $E(Q) = Q\setminus Q'$. Then we clearly have
\begin{equation*}
  |E(Q)| = \left|Q\setminus Q'\right| \geq \left(1-\frac{1}{2^n}\right)\left|Q\right|.
\end{equation*}
Since the supports of $f$ and $g$ do not intersect $E(Q)$, we even have that
\begin{equation*}
    \langle f \rangle_{r,Q'} = 2^{n/r}\langle f \rangle_{r,Q} \quad \mbox{and} \quad \langle g \rangle_{s',Q'} = 2^{n/s'}\langle g \rangle_{s',Q},
\end{equation*}
and
\begin{equation}\label{stoppingtimenegativep}
    \langle f \rangle_{r,Q''} \leq 2^{n/r}\langle f \rangle_{r,Q} \quad \mbox{and} \quad \langle g \rangle_{s',Q''} \leq 2^{n/s'}\langle g \rangle_{s',Q},
\end{equation}
whenever the dyadic cube $Q'' \subseteq Q$ is not a proper subset of $Q'$.
Furthermore, it will be useful to denote by $\tilde{\mathcal{S}}(Q)$ the collection of all dyadic sub-cubes of $Q'$ and $\mathcal{S}(Q) = \{Q'\}$. Thus, $\mathcal{S}(Q)$ is the (singleton) collection of maximal cubes of $\tilde{\mathcal{S}}(Q)$ with respect to inclusion. Although defined differently here for $p<0$ from the case $p\geq0$, the two variants of $\tilde{\mathcal{S}}(Q)$ and $\mathcal{S}(Q)$ play the same role in the proof which follows.

Defining $\mathcal{S}^\omega = \cup_{p\in\Z} \mathcal{S}^\omega_p$, we see that if two cubes $Q$ and $Q'$ from $\mathcal{S}^{\omega}$ intersect, then, as dyadic cubes, one is contained in the other, say $Q' \subseteq Q$. It is clear from the construction that if $Q' \subseteq Q$ then $E(Q') \subseteq Q' \subseteq Q \setminus E(Q)$, and so $E(Q) \cap E(Q') = \emptyset$. Thus we see that $\mathcal{S}^{\omega}$ is $\frac{1}{2}$-sparse. As, by Theorem~1.3 in \cite{Hanninen}, a finite collection of sparse sets is sparse, then even
\begin{equation*}
    \mathcal{S} := \bigcup_{\omega \in \{-1,0,1,\}^n}\mathcal{S}^{\omega}
\end{equation*}
is sparse.

Consider again a fixed $\omega$. The complement of the collection $\tilde{\mathcal{S}}(Q)$ within the collection of all dyadic sub-cubes of a cube $Q$ is
\begin{equation*}
    \mathcal{E}(Q) := \left\{Q' \in \mathcal{D}^\omega \colon Q' \subseteq Q, Q' \not\in \tilde{\mathcal{S}}(Q) \right\},
\end{equation*}
so
\begin{equation*}
    \mathcal{E}(Q) \cup \tilde{\mathcal{S}}(Q)
\end{equation*}
is a partition of all the dyadic sub-cubes of $Q$. Note that since $Q \not\in \mathcal{S}(Q)$, we always have $Q \in \mathcal{E}(Q)$. We now want to show that
\begin{equation} \label{inclusion}
    \bigcup_{p \in \Z} \bigcup_{Q \in \mathcal{S}^\omega_p} \mathcal{E}(Q) \cup \tilde{\mathcal{S}}(Q)
    =  \bigcup_{p \in \Z} \bigcup_{Q \in \mathcal{S}^\omega_p} \mathcal{E}(Q).
\end{equation}
To do this we assume that a cube $Q'$ belongs to $\tilde{\mathcal{S}}(Q)$ with $Q \in \mathcal{S}^\omega_p$. Then either $Q'$ is maximal in $\tilde{\mathcal{S}}(Q)$ (that is, $Q' \in \mathcal{S}(Q)$) or not. If $Q'$ is maximal then $Q' \in \mathcal{S}^\omega_{p+1}$ and $Q' \in \mathcal{E}(Q')$. If $Q'$ is not maximal then there exists a $Q'' \in \mathcal{S}^\omega_{p+1}$ such that $Q \supseteq Q'' \supsetneq Q'$, and either $Q' \in \mathcal{E}(Q'')$ or $Q' \in \tilde{\mathcal{S}}(Q'')$. But it is only in this later case that we have not obtained an outcome we desire to prove \eqref{inclusion}. However, this later case, where $Q'' \in \mathcal{S}^\omega_{p+1}$ and $Q' \in \tilde{\mathcal{S}}(Q'')$, is exactly our original assumption with $Q$ replaced by the strictly smaller cube $Q''$ (and $p$ replaced by $p+1$). Since there are only a finite number of dyadic cubes $Q''$ such that $Q' \subsetneq Q'' \subseteq Q$, we only need repeat this argument a finite number of times to find a cube $Q''$ and an index $q$ such that $Q' \in \mathcal{E}(Q'')$ with $Q'' \in \mathcal{S}^\omega_q$. Consequently, \eqref{inclusion} is proved.

Since the union of cubes of side length $2^{-k}$ in \eqref{inclusion} contains the support of $f$, we can write
\begin{equation}\label{formfunctiondecomp}
    f = \sum_{\omega \in \{-1,0,1,\}^n} \sum_{Q\in\mathcal{S}^\omega} \sum_{Q'_k\in\mathcal{E}(Q)} f\chi_{\frac{1}{3}Q'_k}
\end{equation}
for each $k\in\Z$, where the innermost sum is only taken over $Q'_k\in\mathcal{E}(Q)$ which have side length $2^{-k}$.

With the sparse collection $\mathcal{S}$ and decomposition \eqref{formfunctiondecomp} in hand, we can proceed to the sparse form bound. Recall from the statement of the proposition that the function $\iota$ connects $T_j$ to cubes of size $2^{-\iota(j)}$. Using \eqref{formfunctiondecomp} with $k=\iota(j)$, Hölder's inequality, and both the support properties of $T_j$ and \eqref{intro:Lr_Ls_averagebounds}, we can write
\begin{align*}
     |\langle T(f), g\rangle| &\leq \sum_{j=1}^\infty  \left|\langle T^{j}(f), g \rangle\right| \\
    &=  \sum_{j=1}^\infty \left|\left\langle T^{j}\left(\sum_{\omega \in \{-1,0,1,\}^n} \sum_{Q\in\mathcal{S}^\omega} \sum_{Q'_{\iota(j)}\in\mathcal{E}(Q)} f\chi_{\frac{1}{3}Q'_{\iota(j)}}\right), g\right\rangle\right| \\
    &\leq  \sum_{j=1}^\infty\sum_{\omega \in \{-1,0,1,\}^n} \sum_{Q\in\mathcal{S}^\omega} \sum_{Q'_{\iota(j)}\in\mathcal{E}(Q)} \left\langle \left|T^{j}\left( f\chi_{\frac{1}{3}Q'_{\iota(j)}}\right)\right|, \left|g\right| \right\rangle \\
    &\leq  \sum_{j=1}^\infty\sum_{\omega \in \{-1,0,1,\}^n} \sum_{Q\in\mathcal{S}^\omega} \sum_{Q'_{\iota(j)}\in\mathcal{E}(Q)} \left|Q'_{\iota(j)}\right|\left\langle T^{j}\left(f\chi_{\frac{1}{3}Q'_{\iota(j)}}\right)\right\rangle_{s,Q'_{\iota(j)}} \left\langle g\right\rangle_{s',Q'_{\iota(j)}} \\
    &\leq  \sum_{j=1}^\infty\sum_{\omega \in \{-1,0,1,\}^n} \sum_{Q\in\mathcal{S}^\omega} \sum_{Q'_{\iota(j)}\in\mathcal{E}(Q)} B_{j} \left|Q'_{\iota(j)}\right|\left\langle f \right\rangle_{r,Q'_{\iota(j)}} \left\langle g\right\rangle_{s',Q'_{\iota(j)}}.
\end{align*}
But for each $Q'_{k} \in \mathcal{E}(Q)$ with $Q \in \mathcal{S}^{\omega}$, we know $\langle f \rangle_{r,Q'_{k}} \leq \max\{4,2^n\}^{1/r}\langle f \rangle_{r,Q}$ and $\langle g \rangle_{s',Q'_{k}} \leq \max\{4,2^n\}^{1/s'}\langle g \rangle_{s',Q}$. Indeed, this follows for $Q \in \mathcal{S}^\omega_p$ by the maximality of cubes in $\mathcal{S}(Q)$ in the case $p\geq0$, and by \eqref{stoppingtimenegativep} with $Q'' = Q'_k$ in the case $p<0$. Moreover, for fixed $k$, the collection of all $Q'_{k} \in \mathcal{E}(Q)$ is a pairwise disjoint collection of subsets of $Q$, so (with $k=\iota(j)$) we have
\begin{align*}
    & \sum_{Q'_{\iota(j)}\in\mathcal{E}(Q)} B_{j} \left|Q'_{\iota(j)}\right|\left\langle f \right\rangle_{r,Q'_{\iota(j)}} \left\langle g\right\rangle_{s',Q'_{\iota(j)}} \\
    &\leq \sum_{Q'_{\iota(j)}\in\mathcal{E}(Q)} \max\{4,2^n\}^{\frac{1}{r} + \frac{1}{s'}}B_{j} \left|Q'_{\iota(j)}\right|\left\langle f \right\rangle_{r,Q} \left\langle g\right\rangle_{s',Q} \\
    &\leq \max\{4,2^n\}^{\frac{1}{r} + \frac{1}{s'}}B_{j} \left(\sum_{Q'_{\iota(j)}\in\mathcal{E}(Q)} \left|Q'_{\iota(j)}\right|\right)\left\langle f \right\rangle_{r,Q} \left\langle g\right\rangle_{s',Q} \\
    &\leq \max\{4,2^n\}^{\frac{1}{r} + \frac{1}{s'}}B_{j} \left|Q\right|\left\langle f \right\rangle_{r,Q} \left\langle g\right\rangle_{s',Q}.   
\end{align*}
It follows, then, that
\begin{align*}
    \left|\langle T(f), g\rangle\right| &\leq  \sum_{j=1}^\infty\sum_{\omega \in \{-1,0,1,\}^n} \sum_{Q\in\mathcal{S}^\omega} \sum_{Q'_{\iota(j)}\in\mathcal{E}(Q)} B_{j} \left|Q'_{\iota(j)}\right|\left\langle f \right\rangle_{r,Q'_{\iota(j)}} \left\langle g\right\rangle_{s',Q'_{\iota(j)}} \\
    &\leq  \sum_{j=1}^\infty\sum_{\omega \in \{-1,0,1,\}^n} \sum_{Q\in\mathcal{S}^\omega} \max\{4,2^n\}^{\frac{1}{r} + \frac{1}{s'}}B_{j} \left|Q\right|\left\langle f \right\rangle_{r,Q} \left\langle g\right\rangle_{s',Q} \\
    &= \max\{4,2^n\}^{\frac{1}{r} + \frac{1}{s'}}\left(\sum_{j=1}^{\infty} B_j \right) \sum_{\omega \in \{-1,0,1,\}^n}\sum_{Q \in \mathcal{S}^{\omega}} \left|Q\right|\left\langle f \right\rangle_{r,Q} \left\langle g\right\rangle_{s',Q} \\
    &\lesssim \sum_{Q \in \mathcal{S}}\left|Q\right|\left\langle f \right\rangle_{r,Q} \left\langle g\right\rangle_{s',Q}
\end{align*}
and the proposition is proved.
\end{proof}

A similar method yields the pointwise sparse bound of Proposition~\ref{intro:pwsparse}.

\begin{proof}[Proof of Proposition~\ref{intro:pwsparse}]
We make use of the same sparse collection $\mathcal{S}$ as in the proof of Proposition~\ref{intro:sparseform}, but, as we do not need to involve the function $g$, we can set $g=0$ (so $\tilde{\mathcal{S}}_g(Q) = \emptyset$). Using \eqref{formfunctiondecomp} with $k=\iota(j)$, the sublinearity of $T_j$ and  \eqref{intro:pwparseaveragebound}, we can write
\begin{align*}
     |T(f)(x)| &\leq \sum_{j=1}^\infty \left| T^{j}(f)(x)\right| \\
    &=  \sum_{j=1}^\infty \left|T^{j}\left(\sum_{\omega \in \{-1,0,1,\}^n} \sum_{Q\in\mathcal{S}^\omega} \sum_{Q'_{\iota(j)}\in\mathcal{E}(Q)} f\chi_{\frac{1}{3}Q'_{\iota(j)}}\right)(x)\right| \\
    &\leq \sum_{j=1}^\infty\sum_{\omega \in \{-1,0,1,\}^n} \sum_{Q\in\mathcal{S}^\omega} \sum_{Q'_{\iota(j)}\in\mathcal{E}(Q)} \left|T^{j}\left( f\chi_{\frac{1}{3}Q'_{\iota(j)}}\right)(x)\right| \\
    &\leq \sum_{j=1}^\infty\sum_{\omega \in \{-1,0,1,\}^n} \sum_{Q\in\mathcal{S}^\omega} \sum_{Q'_{\iota(j)}\in\mathcal{E}(Q)} B_{j} \left\langle f \right\rangle_{r,Q'_{\iota(j)}}\chi_{Q'_{\iota(j)}}(x).
\end{align*}
Once again, for each $Q'_{k} \in \mathcal{E}(Q)$ with $Q \in \mathcal{S}^{\omega}$, we know $\langle f \rangle_{r,Q'_{k}} \leq \max\{4,2^n\}^{1/r}\langle f \rangle_{r,Q}$ and, for fixed $k$, the collection of all $Q'_{k} \in \mathcal{E}(Q)$ is a pairwise disjoint collection of subsets of $Q$, so (with $k=\iota(j)$) we have
\begin{align*}
    & \sum_{Q'_{\iota(j)}\in\mathcal{E}(Q)} B_{j} \left\langle f \right\rangle_{r,Q'_{\iota(j)}}\chi_{Q'_{\iota(j)}}(x)\\
    &\leq \sum_{Q'_{\iota(j)}\in\mathcal{E}(Q)} B_{j}\max\{4,2^n\}^{\frac{1}{r}} \left\langle f \right\rangle_{r,Q} \chi_{Q'_{\iota(j)}}(x) \\
    &\leq B_{j}\max\{4,2^n\}^{\frac{1}{r}} \left\langle f \right\rangle_{r,Q}\left(\sum_{Q'_{\iota(j)}\in\mathcal{E}(Q)} \chi_{Q'_{\iota(j)}}(x) \right) \\
    &\leq B_{j}\max\{4,2^n\}^{\frac{1}{r}} \left\langle f \right\rangle_{r,Q} \chi_{Q}(x).
\end{align*}
It follows, then, that
\begin{align*}
    \left|T(f)(x)\right| &\leq \sum_{j=1}^\infty\sum_{\omega \in \{-1,0,1,\}^n} \sum_{Q\in\mathcal{S}^\omega} \sum_{Q'_{\iota(j)}\in\mathcal{E}(Q)} B_{j} \left\langle f \right\rangle_{r,Q'_{\iota(j)}}\chi_{Q'_{\iota(j)}}(x) \\
    &\leq  \sum_{j=1}^\infty\sum_{\omega \in \{-1,0,1,\}^n} \sum_{Q\in\mathcal{S}^\omega} B_{j}\max\{4,2^n\}^{\frac{1}{r}} \left\langle f \right\rangle_{r,Q} \chi_{Q}(x) \\
    &= \max\{4,2^n\}^{\frac{1}{r}}\left(\sum_{j=1}^{\infty} B_j \right) \sum_{\omega \in \{-1,0,1,\}^n} \sum_{Q\in\mathcal{S}^\omega} \left\langle f \right\rangle_{r,Q} \chi_{Q}(x) \\
    &\lesssim \sum_{Q \in \mathcal{S}}\left\langle f \right\rangle_{r,Q} \chi_{Q}(x)
\end{align*}
and the proposition is proved.
\end{proof}

Let us explain how the ideas of Proposition~\ref{intro:sparseform} can be translated to the end-point case. Observe that Proposition~\ref{intro:sparseform} requires two properties of the pieces $T_j$ of the operator $T$: Firstly, that the $T_j$ are bounded operators from $L^r$ to $L^s$ with norms that are appropriately summable in $j$ (that is, \eqref{intro:Lr_Ls_averagebounds} holds); and secondly, that each $T_j$ is a local operator on a scale of the cubes $Q_{\iota(j)}$ with side lengths $2^{-\iota(j)}$. To adapt these ideas to the end-point case, we cannot directly sum estimates made on each $T_j$. Instead, we make use of an induction argument, where the induction is on the size of the cube in our sparse set (in an analogous way to \cite{LMSvideos}). This means some sort of ordering on our sparse set is needed, through which we can perform an induction argument. We collect the notions we need in the following definition.
\begin{definition}\label{def:gradedposet}
    A \emph{partial order} $\preceq$ is a binary relation on certain pairs of elements of a set $\mathcal{S}$ which is \emph{reflexive} ($Q \preceq Q$ for all $Q\in\mathcal{S}$), \emph{antisymmetric} (if $Q' \preceq Q$ and $Q \preceq Q'$ then $Q'=Q$), and \emph{transitive} (if $Q' \preceq Q$ and $Q \preceq Q''$ then $Q' \preceq Q''$). The set $\mathcal{S}$ together with the partial order $\preceq$ is called a \emph{partially ordered set}.

    We write $Q' \prec Q$ when $Q' \preceq Q$ but $Q' \neq Q$, and say $Q$ \emph{covers} $Q'$ and write $Q' \preceqdot Q$ to mean $Q' \prec Q$ but there does not exist a $Q''\in\mathcal{S}$ such that $Q' \prec Q'' \prec Q$.

    For a partially ordered set $\mathcal{S}$ with a partial order $\preceq$, we can define a second partial order $\preceq_{\text{op}}$ on $\mathcal{S}$ called the \emph{opposite partial order} by defining
    \begin{equation*}
        Q \preceq_{\text{op}} Q' \quad\mbox{when} \quad Q' \preceq Q.
    \end{equation*}
    
     A \emph{rank function} $\varrho$ is a function from $\mathcal{S}$ to non-negative integers that is \emph{compatible} with the ordering (if $Q' \prec Q$ then $\varrho(Q') < \varrho(Q)$) and \emph{consistent} with the covering relation (if $Q' \preceqdot Q$, then $\varrho(Q) = \varrho(Q') + 1$). A partially ordered set with a rank function is called a \emph{graded partially ordered set}.
\end{definition}
With these definitions in hand, we formulate two assumptions which will be used in Proposition~\ref{intro:end-pointsparsebdd2}.
\begin{enumerate}[label=(\alph*), resume]
    \item \label{assump:inductionaxiom} The set $\mathcal{S}$ is a partially ordered set of cubes with the partial order $\preceq$. The ordering $Q' \preceq Q$ is defined to hold between two cubes $Q'$ and $Q$ when the inclusion
    \begin{equation} \label{finclusion}
        {\textstyle\frac{1}{3}Q'} \subseteq {\textstyle\frac{1}{3}Q},
    \end{equation}
    holds. Furthermore, for
    \begin{equation} \label{def:f}
        E(Q) := {\textstyle\frac{1}{3}Q} \setminus \bigcup_{Q' \preceqdot Q} {\textstyle\frac{1}{3}Q'},
    \end{equation}
    where the union is assumed to be disjoint and is taken over all $Q'$ which are covered by $Q$, the collection $\{E(Q)\colon Q\in\mathcal{S}\}$ forms a pairwise disjoint collection of sets for which $|E(Q)| \geq \eta|Q|$ for some $0 < \eta \leq 1$, and as such $\mathcal{S}$ is sparse.

    \item \label{assump:flag} We have a rank function $\varrho$ which makes $\mathcal{S}$ into a graded partially ordered set with respect to the opposite partial order to $\preceq$. The rank function $\varrho$ is bounded below by zero and every cube $Q\in\mathcal{S}$ is comparable to some cube $Q_0$ of rank zero. Consequently, for each $Q\in\mathcal{S}$ with $\varrho(Q)=N$, there must exist a sequence of cubes $Q_N, Q_{N-1},\dots,Q_1, Q_0$ such that
    \begin{equation*}
        Q = Q_N \preceqdot Q_{N-1} \preceqdot \dots \preceqdot Q_1 \preceqdot Q_0
    \end{equation*}
    and $\varrho(Q_p) = p$ for $p=0,\dots,N$.

\end{enumerate}    
The $L^r$ to $L^s$ bound \eqref{intro:Lr_Ls_averagebounds} on the pieces of the operator can be replaced by the following $L^r$ to $L^s$ bound for the whole operator between appropriate cubes in the sparse collection.
\begin{enumerate}[label=(\alph*),resume]
    \item \label{assump:lrls} If $Q' \preceqdot Q$, then there exists a constant $A_1$ such that
    \begin{equation*}
        \displaystyle\left\langle T(f\chi_{Q\setminus Q'}) - \left( T(f\chi_{Q\setminus Q'})\right)_{\frac{1}{3}Q'} \right\rangle_{s,\frac{1}{3}Q'} \leq A_1 \left\langle f \right\rangle_{r,Q}.
    \end{equation*}
\end{enumerate}
The induction argument avoids the need for a localisation property on multiple sizes of cubes, but it is useful to assume the following localisation at some level in order to perform the initial step in the induction argument.
\begin{enumerate}[label=(\alph*), resume]
    \item \label{assump:support} The collection $\{\frac{1}{3}Q\}_{\varrho(Q) = 0}$ of dilates of cubes of minimal rank is a countable disjoint collection of dyadic cubes (of the form \eqref{def:dyadic_cube} below for some $k$ and $\omega$) all with the same side length $2^{-k} \geq \ell$,
    \begin{equation*}
        \supp(f) \cup \supp(g) \subseteq \bigcup_{\varrho(Q) = 0} {\textstyle\frac{1}{3}Q},
    \end{equation*}
    where the union is taken over all $Q$ of rank zero, and, again when $\varrho(Q) = 0$,
    \begin{equation*}\label{supportinclusion}
        \supp(T(f\chi_{\frac{1}{3}Q})) \subseteq Q.
    \end{equation*}

\end{enumerate}
We formulate the end result as the following proposition and note that the constant $C$ contained in it does not depend on $\ell$ from assumption \ref{assump:support}.
\begin{proposition}\label{intro:end-pointsparsebdd2}
    Consider $1 \leq r \leq s \leq \infty$. For a linear operator $T$ and each pair of bounded and compactly supported functions $f$ and $g$, assume there exists a collection of cubes $\mathcal{S}$ such that assumptions \ref{assump:inductionaxiom}, \ref{assump:flag}, \ref{assump:lrls} and \ref{assump:support} above hold. Furthermore, define
    \begin{equation}\label{def:gtildeq}
        \tilde{g}_Q(x) = \begin{cases}
            g(x), & \mbox{for $x \in E(Q)$,} \\
            ( g )_{\frac{1}{3}Q'} & \mbox{for $x\in \frac{1}{3}Q'$ when $Q' \preceqdot Q$,}
        \end{cases}
    \end{equation}
    where $(g)_{\frac{1}{3}Q'} = \frac{1}{|{\frac{1}{3}Q'}|} \int_{\frac{1}{3}Q'} g$ and $E(Q)$ is defined in \eqref{def:f}, and assume that there exist constants $A_2$, $A_3$ and $A_4$ such that
    \begin{align}
        \left\langle T(f\chi_{Q})\right\rangle_{r,Q} &\leq A_2 \left\langle f \right\rangle_{r,Q},\label{end-point1} \\
        \left\langle \tilde{g}_Q \right\rangle_{r',Q} &\leq A_3 \left\langle g \right\rangle_{s',Q}\label{end-point2}
        \quad \mbox{and} \\
        \left\langle g \right\rangle_{s',Q'} &\leq A_4 \left\langle g \right\rangle_{s',Q},\label{end-point3}
    \end{align}
    for all $Q\in\mathcal{S}$ and $Q' \preceqdot Q$. Then there exists a constant $C=C(A_1,A_2,A_3,A_4)>0$, such that
    \begin{equation*}
    \left\langle T(f), g\right\rangle \leq C \sum_{Q \in \mathcal{S}} \langle f \rangle_{r,Q}\langle g \rangle_{s',Q}|Q|.
   \end{equation*}
\end{proposition}

We note that \eqref{end-point1}--\eqref{end-point3} are very natural conditions. For example, if $T$ is a bounded operator from $L^r$ to $L^r$, then \eqref{end-point1} holds for all cubes $Q$. We can also easily ensure that \eqref{end-point2} and \eqref{end-point3} hold by constructing $\mathcal{S}$ using a standard Whitney decomposition (see the proof of Lemma~\ref{off-diagonal_extremes} below). 

%%%%%%%%%%%%%%%%%%%%%%%
%%%%%%%%%%%%%%%%%%%%%%%

\begin{proof}[Proof of Proposition~\ref{intro:end-pointsparsebdd2}] 
First, from assumption \ref{assump:support}, we know we can write
\begin{equation*}
    f= \sum_{\varrho(Q) = 0} f\chi_{\frac{1}{3}Q}, \quad g= \sum_{\varrho(Q) = 0} g\chi_{\frac{1}{3}Q} \quad \mbox{and} \quad \supp(T(f\chi_{\frac{1}{3}Q}) \subseteq Q \quad \mbox{when $\varrho(Q) = 0$}.
\end{equation*}
Consequently, denoting $\dist(\frac{1}{3}P,\frac{1}{3}Q) = \inf\{|x-y|\colon x\in \frac{1}{3}P, y\in \frac{1}{3}Q\}$, we can write
\begin{align*}
    \left|\big\langle T(f),g \big\rangle\right| &\leq \sum_{\varrho(Q)=0}  \left|\left\langle T\left(\sum_{\varrho(P)=0}f\chi_{\frac{1}{3}P}\right),g\chi_{\frac{1}{3}Q}\right\rangle\right| \\
    &\leq \sum_{\varrho(Q)=0} \left(\left|\left\langle T\left(\sum_{\substack{\varrho(P)=0\\
                  \dist(\frac{1}{3}P,\frac{1}{3}Q)=0}} f\chi_{\frac{1}{3}P}\right),g\chi_{\frac{1}{3}Q}\right\rangle\right| + \sum_{\substack{\varrho(P)=0\\
                  \dist(\frac{1}{3}P,\frac{1}{3}Q)>0}} \left|\big\langle T(f\chi_{\frac{1}{3}P}),g\chi_{\frac{1}{3}Q}\rangle\right|\right)
\end{align*}
Moreover, since $\{\frac{1}{3}P\}_{\varrho(P)=0}$ are disjoint dyadic cubes with common side length, we can write
\begin{equation*}
    \left\langle T\left(\sum_{\substack{\varrho(P)=0\\
                  \dist(\frac{1}{3}P,\frac{1}{3}Q)=0}} f\chi_{\frac{1}{3}P}\right),g\chi_{\frac{1}{3}Q}\right\rangle
    = \big\langle T(f\chi_{Q}),g\chi_{\frac{1}{3}Q}\big\rangle
\end{equation*}
for each $Q$ with $\varrho(Q)=0$, and, by the support properties of $T(f\chi_{\frac{1}{3}P})$,
\begin{equation*}
    \sum_{\substack{\varrho(P)=0\\
                  \dist(\frac{1}{3}P,\frac{1}{3}Q)>0}} \left|\big\langle T(f\chi_{\frac{1}{3}P}),g\chi_{\frac{1}{3}Q}\big\rangle\right| = \sum_{\substack{\varrho(P)=0\\
                  \dist(\frac{1}{3}P,\frac{1}{3}Q)>0}} \left|\big\langle T(f\chi_{\frac{1}{3}P})\chi_{P},g\chi_{\frac{1}{3}Q}\big\rangle\right|
    = 0,
\end{equation*}
since when $\dist(\frac{1}{3}P,\frac{1}{3}Q)>0$, $\dist(\frac{1}{3}P,\frac{1}{3}Q)$ is at least a third the side length of $Q$, and so $P \cap \frac{1}{3}Q = \emptyset$. Putting this together, we obtain
\begin{equation}\label{basecase}
    \left|\big\langle T(f),g \big\rangle\right| \leq \sum_{\varrho(Q)=0} \left|\big\langle T(f\chi_{Q}),g\chi_{\frac{1}{3}Q}\big\rangle\right|.
\end{equation}

Now fix a $Q \in \mathcal{S}$ of arbitrary rank $p$ and recall the definition of $E(Q)$ in \eqref{def:f}. We can write
    \begin{align*}
        \left|\big\langle T(f\chi_{Q}),g\chi_{\frac{1}{3}Q}\big\rangle\right|
        &\leq  \left|\big\langle T(f\chi_{Q}),\tilde{g}_{Q}\big\rangle\right| + \sum_{Q' \preceqdot Q} \left|\big\langle T(f\chi_{Q}),(g - g_{\frac{1}{3}Q'})\chi_{\frac{1}{3}Q'}\big\rangle\right| \\
        &\leq  \left|\big\langle T(f\chi_{Q}),\tilde{g}_{Q}\big\rangle\right| + \sum_{Q' \preceqdot Q} \left|\big\langle T(f\chi_{Q\setminus Q'}),(g - g_{\frac{1}{3}Q'})\chi_{\frac{1}{3}Q'}\big\rangle\right| \\
        &\quad+ \sum_{Q' \preceqdot Q} \left|\big\langle T(f\chi_{Q'}),(g - g_{\frac{1}{3}Q'})\chi_{\frac{1}{3}Q'}\big\rangle\right|,
    \end{align*}
Since $\mathcal{S}$ is sparse and using \eqref{end-point1} and \eqref{end-point2}, we can estimate
    \begin{align*}
        \left|\big\langle T(f\chi_{Q}),\tilde{g}_{Q}\big\rangle\right| &\leq \big\langle T(f\chi_{Q})\big\rangle_{r,Q} \big\langle \tilde{g}_{Q} \big\rangle_{r',Q}|Q| \\
        &\leq A_2A_3\big\langle f\big\rangle_{r,Q} \big\langle g\big\rangle_{s',Q}|Q|
    \end{align*}
and, by assumption~\ref{assump:lrls}, inequality \eqref{end-point3} and the sparsity of $\mathcal{S}$,
    \begin{align*}
        &\sum_{Q' \preceqdot Q} \left|\big\langle T(f\chi_{Q\setminus Q'}),(g - g_{\frac{1}{3}Q'})\chi_{\frac{1}{3}Q'}\big\rangle\right| \\
        &= \sum_{Q' \preceqdot Q} \left|\big\langle T(f\chi_{Q\setminus Q'}) - ( T(f\chi_{Q\setminus Q'}) )_{\frac{1}{3}Q'} ,(g - g_{\frac{1}{3}Q'})\chi_{\frac{1}{3}Q'}\big\rangle\right| \\
        &\lesssim \sum_{Q' \preceqdot Q} \big\langle T(f\chi_{Q\setminus Q'}) - ( T(f\chi_{Q\setminus Q'}) )_{\frac{1}{3}Q'}\big\rangle_{s,\frac{1}{3}Q'} \big\langle g \big\rangle_{s',{\frac{1}{3}Q'}}|Q'| \\
        &\leq \sum_{Q' \preceqdot Q} A_1\big\langle f\big\rangle_{r,Q} A_4\big\langle g \big\rangle_{s',Q} |Q'| \\
        &\lesssim A_1A_4 \big\langle f\big\rangle_{r,Q} \big\langle g \big\rangle_{s',Q}\sum_{Q' \preceqdot Q} |E(Q')| \\
        &\leq A_1A_4 \big\langle f\big\rangle_{r,Q} \big\langle g \big\rangle_{s',Q} |Q|,
    \end{align*}
so together we have
\begin{equation*}
    \left|\big\langle T(f\chi_{Q}),g\chi_{\frac{1}{3}Q}\big\rangle\right| \leq C(A_2A_3 + A_1A_4)\big\langle f\big\rangle_{r,Q}\big\langle g \big\rangle_{s',Q}|Q| + \sum_{Q' \preceqdot Q} \left|\big\langle T(f\chi_{Q'}),g\chi_{\frac{1}{3}Q'}\big\rangle\right|
\end{equation*}
for each $Q \in \mathcal{S}$. Summing the last inequality over all $Q$ of rank $p$ we obtain
\begin{equation}
\begin{aligned}\label{inductivestep}
    &\sum_{\varrho(Q)=p}\left|\big\langle T(f\chi_{Q}),g\chi_{\frac{1}{3}Q}\big\rangle\right| \\
    &\leq C_0\sum_{\varrho(Q)=p}\big\langle f\big\rangle_{r,Q}\big\langle g \big\rangle_{s',Q}|Q| + \sum_{\varrho(Q)=p}\sum_{Q' \preceqdot Q} \left|\big\langle T(f\chi_{Q'}),g\chi_{\frac{1}{3}Q'}\big\rangle\right| \\
    &\leq C_0\sum_{\varrho(Q)=p}\big\langle f\big\rangle_{r,Q}\big\langle g \big\rangle_{s',Q}|Q| + \sum_{\varrho(Q')=p+1} \left|\big\langle T(f\chi_{Q'}),g\chi_{\frac{1}{3}Q'}\big\rangle\right|,
\end{aligned}
\end{equation}
where we wrote $C_0 = C(A_2A_3 + A_1A_4)$. Combining \eqref{basecase} and $N$ applications of \eqref{inductivestep}, we obtain
\begin{equation}
\begin{aligned}\label{repeatedsteps}
    \left|\big\langle T(f),g \big\rangle\right| &\leq \sum_{\varrho(Q)=0} \left|\big\langle T(f\chi_{Q}),g\chi_{\frac{1}{3}Q}\big\rangle\right| \\
    &\leq C_0\sum_{\varrho(Q)=0} \big\langle f\big\rangle_{r,Q}\big\langle g \big\rangle_{s',Q}|Q| + \sum_{\varrho(Q)=1} \left|\big\langle T(f\chi_{Q}),g\chi_{\frac{1}{3}Q}\big\rangle\right| \\
    &\leq C_0\sum_{\varrho(Q)=0,1} \big\langle f\big\rangle_{r,Q}\big\langle g \big\rangle_{s',Q}|Q| + \sum_{\varrho(Q)=2} \left|\big\langle T(f\chi_{Q}),g\chi_{\frac{1}{3}Q}\big\rangle\right| \\
    &\leq C_0\sum_{\varrho(Q)=0,1,2} \big\langle f\big\rangle_{r,Q}\big\langle g \big\rangle_{s',Q}|Q| + \sum_{\varrho(Q)=3} \left|\big\langle T(f\chi_{Q}),g\chi_{\frac{1}{3}Q}\big\rangle\right| \\
    &\,\,\vdots \\
    &\leq C_0\sum_{\varrho(Q)=1,\dots,N-1} \big\langle f\big\rangle_{r,Q}\big\langle g \big\rangle_{s',Q}|Q| + \sum_{\varrho(Q)=N} \left|\big\langle T(f\chi_{Q}),g\chi_{\frac{1}{3}Q}\big\rangle\right|.
\end{aligned}
\end{equation}
Now, by assumption~\ref{assump:inductionaxiom},
\begin{equation*}
    \sum_{Q'\preceqdot Q} |{\textstyle\frac{1}{3}Q'}| = |{\textstyle\frac{1}{3}Q} \setminus E(Q)| \leq (1-\eta)|{\textstyle\frac{1}{3}Q}|,
\end{equation*}
but, by assumption~\ref{assump:flag}, for every $Q'$ of rank $p$ there exists a cube $Q$ of rank $p-1$ such that $Q'\preceqdot Q$, so
\begin{equation*}
    \sum_{\rho(Q')=p} |Q'| \leq \sum_{\substack{Q'\preceqdot Q\\
                  \varrho(Q)=p-1}} |Q'| \leq (1-\eta)\sum_{\rho(Q)=p-1}|Q|.
\end{equation*}
Repeated application of this inequality for $p=1,\dots,N$ gives
\begin{equation}\label{cubestozero}
    \sum_{\rho(Q)=N} |Q| \leq (1-\eta)^N\sum_{\rho(Q)=0}|Q|.
\end{equation}
We can also compute using \eqref{end-point1} that
\begin{align*}
     \left|\big\langle T(f\chi_{Q}),g\chi_{\frac{1}{3}Q}\big\rangle\right| &\lesssim  \big\langle T(f\chi_{Q}) \big\rangle_{r,\frac{1}{3}Q} \big\langle g\chi\big\rangle_{r',\frac{1}{3}Q}|Q| \\
     &\lesssim \big\langle f \big\rangle_{r,Q} \big\langle g\chi\big\rangle_{r',\frac{1}{3}Q}|Q| \\
     &\lesssim \big\langle f \big\rangle_{r,Q} \big\langle g\chi\big\rangle_{r',\frac{1}{3}Q}|Q| \\
     &\lesssim \| f \|_{L^\infty}\| g \|_{L^\infty} |Q|,
\end{align*}
so, using \eqref{cubestozero}, we see that
\begin{equation*}
    \sum_{\varrho(Q)=N} \left|\big\langle T(f\chi_{Q}),g\chi_{\frac{1}{3}Q}\big\rangle\right| \lesssim \| f \|_{L^\infty}\| g \|_{L^\infty} \sum_{\varrho(Q)=N} |Q| \lesssim (1-\eta)^N\| f \|_{L^\infty}\| g \|_{L^\infty} \sum_{\varrho(Q)=0} |Q| \to 0
\end{equation*}
as $N\to\infty$. Since
\begin{equation*}
    \lim_{N\to\infty} \sum_{\varrho(Q)=1,\dots,N-1} \big\langle f\big\rangle_{r,Q}\big\langle g \big\rangle_{s',Q}|Q| = \sum_{Q\in\mathcal{S}} \big\langle f\big\rangle_{r,Q}\big\langle g \big\rangle_{s',Q}|Q|,
\end{equation*}
letting $N\to\infty$ in \eqref{repeatedsteps} proves the proposition.
\end{proof}

%%%%%%%%%%%%%%%%%%%%%%
%%%%%%%%%%%%%%%%%%%%%%
%%%%%%%%%%%%%%%%%%%%%%
\section{Rough and Smooth Pseudodifferential Operators}
\label{sec:3}

In this section we apply Proposition~\ref{intro:pwsparse} to prove Theorems~\ref{thm:rough-symbols} and \ref{thm:rough-amplitudes}. Our methods here are largely the same as those in \cite{BC}, but we include them for sake of presenting a complete argument.

We begin, just as in \cite{BC}, by decomposing the pseudodifferential operator $T_a$, with $a \in L^\infty A^m_\rho$, in both frequency and space. First, for the frequency decomposition, we decompose the operator into Littlewood-Paley pieces. We introduce a smooth cut-off function $\psi_0$ which is supported on the unit ball centred at the origin, and equal to one on the ball of radius $\frac{1}{2}$.\label{page:psi_0} We then define $\psi(\xi) := \psi_0(2^{-1}\xi) - \psi_0(\xi)$ and $\psi_k(\xi) := \psi(2^{1-j}\xi)$ for $j \in \N := \{1,2,3,\dots\}$, so
\begin{equation*}
    \sum_{j=0}^\infty \psi_k(\xi) = 1
\end{equation*}
for all $\xi \in \R^n$. We then decompose the operator as
\begin{equation*}
    T_a = \sum_{j=0}^\infty T_a^j,
\end{equation*}
where
\begin{equation}\label{def:taj}
    T_a^jf(x) = \int_{\R^n}\int_{\R^n} a(x,y,\xi)\psi_j(\xi) f(y) e^{i(x-y)\cdot\xi} dy d\xi = \int_{\R^n} K_a^{j}(x,y,x-y) f(y) dy
\end{equation}
for all $j \in \N_0 := \{0\} \cup \N$, and
\begin{equation*}
    K_a^{j}(x,y,z) = \int_{\R^n} a(x,y,\xi)\psi_j(\xi) e^{iz\cdot\xi} d\xi.
\end{equation*}
Secondly, we localise in the $x$-variable. To do this, we make use of the same functions $\psi$ and $\psi_0$.  We define $T_{a,\nu}^{j,\ell}$ so that for each $j \in\N_0$
\begin{equation*}
    T_a^j = \sum_{\ell=0}^\infty T_{a,\nu}^{j,\ell},
\end{equation*}
where, for given $\nu \in [0,1]$ and $\ell \in \N$,
\begin{equation}
\begin{aligned}\label{def:Tjlanu}
    T_{a,\nu}^{j,\ell}f(x) &:= \int_{\R^n} a(x,y,\xi)\psi_j(\xi) \int_{\R^n} f(y) \psi(2^{j\nu-\ell}(x-y)) e^{i(x-y)\cdot\xi} dy d\xi \\
    &= \int_{\R^n} K_{a,\nu}^{j,\ell}(x,y,x-y) f(y) dy
\end{aligned}
\end{equation}
and
\begin{equation*}
    K_{a,\nu}^{j,\ell}(x,y,z) = \int_{\R^n} a(x,y,\xi)\psi_j(\xi) \psi(2^{j\nu-\ell}z) e^{iz\cdot\xi} d\xi.
\end{equation*}
and for $\ell = 0$ we have
\begin{equation}
\begin{aligned}\label{def:Tj0anu}
    T_{a,\nu}^{j,0}f(x) &:= \int_{\R^n} a(x,y,\xi)\psi_j(\xi) \int_{\R^n} f(y) \psi_0(2^{j\nu-1}(x-y)) e^{i(x-y)\cdot\xi} dy d\xi \\
    &= \int_{\R^n} K_{a,\nu}^{j,0}(x,y,x-y) f(y) dy
\end{aligned}
\end{equation}
and
\begin{equation*}
    K_{a,\nu}^{j,0}(x,y,z) = \int_{\R^n} a(x,y,\xi)\psi_j(\xi) \psi_0(2^{j\nu-1}z) e^{iz\cdot\xi} d\xi.
\end{equation*}
The parameter $\nu$ in this decomposition will be chosen to be slightly less than $\rho$ in order to prove Theorems~\ref{thm:rough-symbols} and \ref{thm:rough-amplitudes}. By leaving it unfixed now, we hope to elucidate how the numerology in the statements of Theorems~\ref{thm:rough-symbols} and \ref{thm:rough-amplitudes} appears.
%%%%%%%%%%%%%%%
\begin{lem}[rough amplitudes, non-local part]\label{amp_lrls}
For $a \in L^\infty A^m_\rho$, with $m, \rho \in \R$ and $\nu < \rho$, we have the estimate
\begin{equation*}
    \|T_{a,\nu}^{j,\ell}\|_{L^r\to L^s} \lesssim 2^{-N_1j}2^{-N_2\ell}2^{-n(\ell - j\nu)(\frac{1}{s} - \frac{1}{r})}
\end{equation*}
for $N_1,N_2>0$, $\ell \geq 1$ and $1 \leq r \leq s \leq \infty$. Furthermore, for $N_1,N_2>0$, $\ell \geq 1$ and any cube $Q$ of side length $2^{\lfloor \ell-j\nu \rfloor}$,
\begin{align*}
    \left|T_{a,\nu}^{j,\ell}(f\chi_{\frac{1}{3}Q})(x)\right| &\lesssim 2^{-N_1j}2^{-N_2\ell} \langle f\rangle_{r,Q}\chi_Q(x) \quad \mbox{and} \\
    \left\langle T_{a,\nu}^{j,\ell}(f\chi_{\frac{1}{3}Q})\right\rangle_{s,Q} &\lesssim 2^{-N_1j}2^{-N_2\ell} \langle f\rangle_{r,Q} \quad \mbox{when $r \leq s$.}
\end{align*}
\end{lem}
%%%%%%%%%%%%%%%
\begin{proof}
We use the identity $\Delta_\xi e^{iz\cdot\xi} = |z|^2e^{iz\cdot\xi}$ to compute
\begin{align*}
    |K_{a,\nu}^{j,\ell}(x,y,z)| &= \left|\int_{\R^n} a(x,y,\xi)\psi_j(\xi) \psi(2^{j\nu-\ell}z) \frac{\Delta_\xi^N e^{iz\cdot\xi}}{|z|^{2N}} d\xi\right| \\
    &= \left|\int_{\R^n} \Delta_\xi^N(a(x,y,\xi)\psi_j(\xi)) \frac{\psi(2^{j\nu-\ell}z)}{|z|^{2N}} e^{iz\cdot\xi} d\xi\right| \\
    &\lesssim 2^{jn}2^{j(m-2\rho N)}2^{2N(j\nu-\ell)}\\
    & = 2^{j(n+m+2N(\nu-\rho))}2^{-2N\ell},
\end{align*}
from the $\xi$ and $z$-support properties of the kernel. Note that we need $\ell > 0$ here. We can also see that $K_a^{j,\ell}(x,z)$ has $z$-support contained in the set $|z| \lesssim 2^{\ell - j\nu}$, so together with the previous estimate, we obtain that
\begin{align*}
    \sup_{x\in \R^n}\left(\int_{\R^n} |K_{a,\nu}^{j,\ell}(x,y,x-y)|^pdy\right)^\frac{1}{p} &\simeq \sup_{y\in \R^n}\left(\int_{\R^n} |K_{a,\nu}^{j,\ell}(x,y,x-y)|^pdx\right)^\frac{1}{p} \\
    &\lesssim 2^{j(n+m+2N(\nu-\rho))}2^{-2N\ell}2^{\frac{n}{p}(\ell - j\nu)}\\
    &= 2^{j(n+m+2N(\nu-\rho))}2^{-2N\ell}2^{n(\ell - j\nu)}2^{-\frac{n}{p'}(\ell - j\nu)}\\ 
    &= 2^{j(n(1-\nu)+m+2N(\nu-\rho))+\ell(n-2N)}2^{-\frac{n}{p'}(\ell - j\nu)}.
\end{align*}
If $\nu < \rho$, we can choose $N \geq \max\bigg\{\frac{n(1-\nu)+m+N_1}{2(\rho-\nu)},\frac{n+N_2}{2}\bigg\}$, so
\begin{equation*}
    j(n(1-\nu)+m+2N(\nu-\rho))+\ell(n-2N) \leq -N_1j - N_2\ell.
\end{equation*}
We can therefore conclude via the Schur test that
\begin{align*}
   \|T_{a,\nu}^{j,\ell}\|_{L^r\to L^s} &\leq \Big(\sup_{x\in \R^n}\int_{\R^n} |K_{a,\nu}^{j,\ell}(x,y,x-y)|^p dy\Big)^\frac{1}{p} + \Big(\sup_{y\in \R^n}\int_{\R^n} |K_{a,\nu}^{j,\ell}(x,y,x-y)|^p dx\Big)^\frac{1}{p} \\
   &\lesssim 2^{-N_1j}2^{-N_2\ell}2^{-n(\ell - j\nu)(\frac{1}{s} - \frac{1}{r})}
\end{align*}
for $\frac{1}{s} + 1 = \frac{1}{p} + \frac{1}{r}$ and $p\geq 1$. This proves the first statement.

The second and third statement follow from the first and the support properties of the kernel $K_{a,\nu}^{j,\ell}(x,y,z)$ in the $z$-variable, which imply that if a function $f$ has support in a cube of side length $2^{\lfloor \ell-j\nu \rfloor}$, then $T_{a,\nu}^{j,\ell}(f)$ is supported in the double of that cube. 
\end{proof}

%%%%%%%%%%%%%%%
\begin{lem}[rough amplitudes, local part]\label{l1linftyamp}
For $a \in L^\infty A^m_\rho$, with $m, \rho \in \R$, we have that the estimates
\begin{equation*}
    \|T_{a,\nu}^{j,0}\|_{L^1\to L^\infty} \lesssim 2^{j(n+m)},
\end{equation*}
and
\begin{equation*}
    \left|T_{a,\nu}^{j,0}(f\chi_{\frac{1}{3}Q})(x)\right| \lesssim 2^{j(n(1-\nu)+m)} \langle f\rangle_{1,Q}\chi_{Q}(x)
\end{equation*}
hold uniformly in $\nu \in \R$, where $Q$ is any cube of side length $2^{\lfloor -j\nu \rfloor}$.
\end{lem}
%%%%%%%%%%%%%%%
\begin{proof}
From the definition of $K_{a,\nu}^{j,0}(x,y,z)$, we have that
\begin{align*}
   |K_{a,\nu}^{j,0}(x,y,z)| &\leq \int_{\R^n} |a(x,y,\xi)\psi_j(\xi)\psi_0(2^{j\nu-1}z) e^{i\xi\cdot z}| d\xi\\
   &\lesssim 2^{jm} 2^{jn} = 2^{j(m+n)}.
\end{align*}
where the last inequality uses the fact that $a(x,y,\xi)\in L^\infty A^m_\rho$ and the measure of the $\xi$-support of $\psi_j(\xi)$ is $2^{jn}$. 

Therefore,
\begin{equation*}\label{l1linftynonlocal}
\begin{aligned}
   \left| T_{a,\nu}^{j,0}f(x)\right| &\leq \int_{\R^n}\left|K_{a,\nu}^{j,0}(x,y,x-y) f(y) \right| dy\\
   &\lesssim 2^{j(m+n)}\|f\|_{L^1},
\end{aligned}
\end{equation*}
and so $\|T_{a,\nu}^{j,0}\|_{L^1\to L^\infty} \lesssim 2^{j(m+n)}$. The second estimate follows from this and the support properties of $K_{a,\nu}^{j,0}$.
\end{proof}

We can now apply Proposition~\ref{intro:pwsparse} to prove Theorem~\ref{thm:rough-amplitudes}. Lemma~\ref{amp_lrls} with $\nu < \rho$, $N_1 \geq -n(1-\nu) - m$ and $N_2=1$, and Lemma~\ref{l1linftyamp} together say that we have the estimate
\begin{equation*}
    \left|T_{a,\nu}^{j,\ell}(f\chi_{\frac{1}{3}Q_{\iota(j,\ell)}})(x)\right| \lesssim 2^{-\ell}2^{j(n(1-\nu)+m)} \langle f\rangle_{1,Q_{\iota(j,\ell)}}\chi_{Q_{\iota(j,\ell)}}(x)
\end{equation*}
for all $(j,\ell) \in \N_0\times\N_0$ and any cube of side length $2^{-\iota(j,\ell)}$ with $\iota(j,\ell) = -\lfloor \ell-j\nu \rfloor$. Choosing $\nu < \rho$ such that $n(1-\nu)+m < 0$, which is possible when $m < -n(1-\rho)$, we can ensure that the constants in this inequality are summable in $j$ and $\ell$. 

Thus we have proved \eqref{intro:pwparseaveragebound} with the countable index set $I = \N_0\times\N_0$ (and with the indices $(j,\ell)$ in place of $j$), and so Proposition~\ref{intro:pwsparse} gives Theorem~\ref{thm:rough-amplitudes}.

%%%%%%%%%%%%%%%
%%%%%%%%%%%%%%
\begin{lem}[rough symbols, local part]\label{symbolstopline}
    For $1 \leq r \leq 2$ and $a \in L^\infty S^m_\rho$, with $m, \rho \in \R$, we have that the estimates
\begin{equation*}
    \left\|T_{a,\nu}^{j,0}\right\|_{L^{r}\to L^\infty} \lesssim 2^{j(m+\frac{n}{r})}
\end{equation*}
and
\begin{equation*}
    \left|T_{a,\nu}^{j,0}(f\chi_{\frac{1}{3}Q})(x)\right| \lesssim 2^{j(\frac{n}{r}(1-\nu)+m)} \langle f\rangle_{r,Q}\chi_Q(x)
\end{equation*}
hold uniformly in $\nu \in \R$, where $Q$ is any cube of side length $2^{\lfloor -j\nu \rfloor}$.
\end{lem}
\begin{proof}
Here our method diverges from \cite{BC}. First, we make use of a $TT^*$-argument with $T=T^j_a$ from \eqref{def:taj}. A calculation shows that
 \begin{equation*}
 T_a^j(T_a^j)^*f(x) = \int_{\R^n}\int_{\R^n} a(x,\xi)\psi_j(\xi)\overline{a(y,\xi)\psi_j(\xi)}  f(y) e^{i(x-y)\cdot\xi} dy d\xi,
 \end{equation*}
 so the composition $T_a^j\circ(T_a^j)^*$ is itself a Littlewood-Paley piece of a pseudodifferential operator with amplitude in $L^\infty A^{2m}_\rho$. Therefore, Lemma~\ref{l1linftyamp} tells us that $\|T_a^j\circ(T_a^j)^*\|_{L^1\to L^\infty} \lesssim 2^{j(2m+n)}$ and, since
\begin{align*}
    \|(T_a^j)^*f\|_{L^2}^2 &= \langle (T_a^j)^*f,(T_a^j)^*f\rangle = \langle T_a^j\circ(T_a^j)^*f,f\rangle \\
    &\leq \|T_a^j\circ(T_a^j)^*f\|_{L^\infty}\|f\|_{L^1} \leq \|T_a^j\circ(T_a^j)^*\|_{L^1 \to L^\infty}\|f\|_{L^1}^2,
\end{align*}
we have that $\|(T_a^j)^*\|_{L^1\to L^2} \lesssim 2^{j(m+\frac{n}{2})}$, so even $\|T_a^j\|_{L^2 \to L^\infty} \lesssim 2^{j(m+\frac{n}{2})}$. This, together with Lemma~\ref{amp_lrls} applied with $N_1 \geq -m-\frac{n}{2}-n\nu(\frac{1}{r}-\frac{1}{s})$ and $N_2 \geq n(\frac{1}{r}-\frac{1}{s})+1$, gives
\begin{align*}
    \|T_{a,\nu}^{j,0}\|_{L^2\to L^\infty} &\leq \|T_a^{j}\|_{L^2\to L^\infty} + \sum_{\ell = 1}^\infty \|T_{a,\nu}^{j,\ell}\|_{L^2\to L^\infty} \\
    &\lesssim 2^{j(m+\frac{n}{2})} + \sum_{\ell = 1}^\infty 2^{j(m+\frac{n}{2})}2^{-\ell} \lesssim 2^{j(m+\frac{n}{2})},
\end{align*}
proving the first inequality of the lemma when $r=2$.

When $r=1$ the lemma is a special case of Lemma~\ref{l1linftyamp}, and for the remaining cases $1 < r < 2$, we can apply the Riesz-Thorin Interpolation Theorem to obtain
\begin{align*}
\left\|T_{a,\nu}^{j,0}\right\|_{L^{r} \to L^\infty} &\lesssim\left\|T_{a,\nu}^{j,0}\right\|_{L^1 \to L^\infty}^{\frac{2}{r}-1} \left\|T_{a,\nu}^{j,0}\right\|_{L^2 \to L^\infty}^{2-\frac{2}{r}}.\\
&\lesssim 2^{j(m+n)(\frac{2}{r}-1)}2^{j(m+\frac{n}{2})(2-\frac{2}{r})} = 2^{j(m+\frac{n}{r})}.
\end{align*}
The second inequality follows directly from the support properties of $K_{a,\nu}^{j,0}$.
\end{proof}

We can now see how Proposition~\ref{intro:pwsparse} proves Theorem~\ref{thm:rough-symbols}. Indeed, once again Lemma~\ref{amp_lrls}, now with $N_1 \geq -\frac{n}{r}(1-\nu) - m$, and again $N_2=1$ and $\nu < \rho$, and this time Lemma~\ref{symbolstopline}, say that we have
\begin{equation*}
    \left|T_{a,\nu}^{j,\ell}(f\chi_{\frac{1}{3}Q_{\iota(j,\ell)}})(x)\right| \lesssim 2^{-\ell}2^{j(\frac{n}{r}(1-\nu)+m)} \langle f\rangle_{r,Q_{\iota(j,\ell)}}\chi_{Q_{\iota(j,\ell)}}(x)
\end{equation*}
for all $(j,\ell) \in \N_0\times\N_0$ and any cube of side length $2^{-\iota(j,\ell)}$ with $\iota(j,\ell) = -\lfloor \ell-j\nu \rfloor$. Now choosing $\nu < \rho$ such that $\frac{n}{r}(1-\nu)+m < 0$, which is possible when $m < -\frac{n}{r}(1-\rho)$, we can ensure that the constants in this inequality are summable in $j$ and $\ell$. Thus we have proved \eqref{intro:pwparseaveragebound} with the countable index set $I = \N_0\times\N_0$ (and with the indices $(j,\ell)$ in place of $j$) and with it Theorem~\ref{thm:rough-symbols}.

%%%%%%%%%%%%%%%%%%%%%%
\begin{lem}[smooth symbols, local part]\label{lrlssmoothsymbols}
    For $1 \leq r \leq 2 \leq r' \leq s$, $m\in \R$, $\rho \in (0,1]$, $\delta \in [0,1)$ and $\nu \in [0,1)$, set $\mu = \max\{0,(\delta-\rho)/s,(\nu-\rho)/s\}$. For $a \in S^m_{\rho,\delta}$, we have the estimates
\begin{equation*}
    \left\|T_{a,\nu}^{j,0}\right\|_{L^{r}\to L^{s}} \lesssim 2^{j(m+n\mu+n(\frac{1}{r}-\frac{1}{s}))},
\end{equation*}
    and
    \begin{equation*}
    \left\langle T_{a,\nu}^{j,0}(f\chi_{\frac{1}{3}Q})\right\rangle_{s,Q} \lesssim 2^{j(m+n\mu+n(1-\mu)(\frac{1}{r}-\frac{1}{s}))}\left\langle f \right\rangle_{r,Q},
    \end{equation*}
    where $Q$ is any cube of side length $2^{\lfloor -j\nu \rfloor}$. Furthermore, the support of $T_{a,\nu}^{j,0}(f\chi_{\frac{1}{3}Q})$ is contained in the cube $Q$.
\end{lem}
\begin{proof}
Set $\delta_0 = \max\{\delta,\nu\}$ and $\lambda_0 = \max\{0,(\delta_0-\rho)/2\}$. The operator $T^{j,0}_{a,\nu}$ is a pseudodifferential operator with symbol in $S^{-n\lambda_0}_{\rho,\delta_0}$ with semi-norms of size $2^{j(m+n\lambda_0)}$. From page~766 in \cite{Hounie} (alternatively, see Theorem~2.2 in \cite{AlvarezHounie90} for the same statement), we know that
\begin{equation}\label{l2boundedness}
     \left\|T_a^{j,0}\right\|_{L^2\to L^2} \lesssim 2^{j(m+n\lambda_0)}.
 \end{equation}
We now wish to interpolate between this $L^2 \to L^2$ boundedness and $L^{r_0} \to L^\infty$ boundedness available from Lemma~\ref{symbolstopline} for some appropriate $r_0$. This requires finding an $r_0 \geq 1$ such that
\begin{equation*}
    \frac{1}{r} = \frac{\left(\frac{2}{s}\right)}{2} + \frac{\left(1-\frac{2}{s}\right)}{r_0},
\end{equation*}
which is possible provided that $s' \leq r$. Applying the Riesz-Thorin interpolation theorem between \eqref{l2boundedness} and Lemma~\ref{symbolstopline}, with $r$ replaced by $r_0$ gives
\begin{align*}
\left\| T_{a,\nu}^{j,0} \right\|_{L^r \to L^s} &\lesssim \left\| T_{a,\nu}^{j,0} \right\|_{L^{r_0} \to L^\infty}^{1-\frac{2}{s}} \left\| T_{a,\nu}^{j,0} \right\|_{L^2 \to L^2}^{\frac{2}{s}}.
\\
&\lesssim 2^{j(m + \frac{n}{r_0})(1-\frac{2}{s})} 2^{j(m+n\lambda_0)\frac{2}{s}} \\
&\lesssim 2^{j(m + n\mu + n(\frac{1}{r} - \frac{1}{s}))}.
\end{align*}
As before, the remaining parts of the lemma follow from the support properties of $K_{a,\nu}^{j,0}$.
\end{proof}

Lemma~\ref{lrlssmoothsymbols} enables us to apply Proposition~\ref{intro:sparseform} and so give an alternative proof of Theorem~1.2 in \cite{BC}. We  apply Lemma~\ref{amp_lrls} with $\nu<\rho$, $N_1 \geq -m-n\mu-n(1-\rho)(\frac{1}{r}-\frac{1}{s})$, where $\mu = \max\{0,(\delta-\rho)/s\}$, and $N_2=1$, and Lemma~\ref{lrlssmoothsymbols} to get that 
\begin{equation*}
    \left\langle T_a^{j,\ell}(f\chi_{\frac{1}{3}Q_{\iota(j,\ell)}})\right\rangle_{s,Q_{\iota(j,\ell)}} \lesssim 2^{-\ell}2^{j(m+n\mu+n(1-\nu)(\frac{1}{r}-\frac{1}{s}))}\left\langle f \right\rangle_{r,Q_{\iota(j,\ell)}},
\end{equation*}
for $1 \leq s' \leq r \leq 2$, all $(j,\ell) \in \N_0\times\N_0$ and any cube of side length $2^{-\iota(j,\ell)}$ with $\iota(j,\ell) = -\lfloor \ell-j\nu \rfloor$. Choosing $\nu < \rho$ such that $m+n\mu+n(1-\rho)(\frac{1}{r}-\frac{1}{s}) < 0$, which is possible when $m < -n\mu-n(1-\rho)(\frac{1}{r}-\frac{1}{s})$, we can ensure that the constants in this inequality are summable in $j$ and $\ell$. Thus we have proved \eqref{intro:Lr_Ls_averagebounds} with the countable index set $I = \N_0\times\N_0$. Since the necessary support properties of $T_a^{j,\ell}(f\chi_{\frac{1}{3}Q_{\iota(j,\ell)}})$ are also given by Lemmas~\ref{amp_lrls} and \ref{lrlssmoothsymbols}, Proposition~\ref{intro:sparseform} shows that when $\rho \in (0,1]$, $\delta \in [0,1)$ and
\begin{equation}\label{deltagrho}
    m < -n(1-\rho)\left(\frac{1}{r} - \frac{1}{s}\right) - \frac{n}{s}\max\{0,\delta-\rho\} \quad \mbox{for} \quad 1\leq s' \leq r\leq 2
\end{equation}
then there exists a constant $C=C(m,\rho,\delta,r,s)$ such that for each pair of bounded functions $f$ and $g$ with compact support, there exists a sparse collection $\mathcal{S}$ such that
\begin{equation*}
    \left|\langle T_a(f),g\rangle\right| \leq C \sum_{Q\in\mathcal{S}} |Q|\langle f \rangle_{r,Q}\langle g \rangle_{s',Q},
\end{equation*}
for $a\in S^m_{\rho,\delta}$. When $\delta < \rho$ the same result holds for the adjoint of $T_a$ and we recover the sparse form bound in the range $1\leq r\leq 2 \leq s \leq \infty$, which is to say, Theorem~1.2 in \cite{BC}.

%%%%%%%%%%%%%%%%%%%%%%
%%%%%%%%%%%%%%%%%%%%%%
%%%%%%%%%%%%%%%%%%%%%%

\section{Sharp Maximal Function Bounds}
\label{sec:4}

In this section we will apply Proposition~\ref{intro:end-pointsparsebdd2} to prove Theorem~\ref{thm:endpoint-pseudos}. In the process of doing this, we improve on previously know pointwise bounds on the composition of a pseudodifferential operator with the sharp maximal function. This was stated as Theorem~\ref{sharpestimatethm}.

We will make use of the following slight generalisation of Lemma 3.2 in \cite{Michalowski_Rule_Staubach_2012}. It is more general from that in \cite{Michalowski_Rule_Staubach_2012} only in the sense that we obtain a localised version \eqref{kernellem2} of the kernel estimate \eqref{kernellem}. Recall that $\psi_0$ is defined on page~\pageref{page:psi_0}.

\begin{lem}\label{kernelestimatelemma}
Let $a\in S^{m}_{\rho, \delta}$, $0 \leq \delta \leq 1$, $0 < \rho \leq 1$ and define $K(x,z)=\int e^{iz\cdot\xi} a(x,\xi) d\xi$ and $\tilde{K}_\ell(x,z)=\psi_0(2^{-\ell}z)\int e^{iz\cdot\xi} a(x,\xi) d\xi$ for $\ell \geq 0$. Then for
\begin{align*}
&\mbox{$|x-x_B|\leq \tau \leq 1$,\quad $\theta \in [0,1]$, \quad $p\in[1,2]$, \quad  $\displaystyle\frac{1}{p}+\frac{1}{p'}=1$,}\\
&\mbox{$\displaystyle m+\frac{n}{p}<h\rho<m+\frac{n}{p}+1$, \quad $\displaystyle\frac{1}{2}< c_1 < 2 c_2 <\infty$ \quad and \quad$j\geq 1$,} 
\end{align*}
we have that
\begin{equation}\label{kernellem}
\left(\int_{A_j}|K(x,x-y)-K(x_B , x_B -y)|^{p'} dy\right)^{\frac{1}{p'}}\lesssim 2^{-jh} \tau^{h(\rho-\theta)-m-\frac{n}{p}}
\end{equation}
and
\begin{equation}\label{kernellem2}
\left(\int_{A_j}|\tilde{K}_\ell(x,x-y)-\tilde{K}_\ell(x_B , x_B -y)|^{p'} dy\right)^{\frac{1}{p'}}\lesssim 2^{-jh} \tau^{h(\rho-\theta)-m-\frac{n}{p}},
\end{equation}
where $A_j = \{y \,|\, c_12^j\tau^\theta \leq |y-x_B| \leq c_22^{j+1}\tau^\theta\}$
\end{lem}
\begin{proof}
The proof of \eqref{kernellem} is that of Lemma 3.2 in \cite{Michalowski_Rule_Staubach_2012}. To prove \eqref{kernellem2}, we can use the same method. Splitting the integral as in \cite{Michalowski_Rule_Staubach_2012} as $I_1$, $I_2$ and $I_3$, the estimate of $I_1$ and $I_2$ is identical, only needing the additional bound $|\psi_0(2^{-\ell}z)| \leq 1$. Estimating $I_3$ is similar, but a third term
\begin{equation*}
    I_{3,3} = \sum_{k\leq k_0} \left( \int_{A_j} \left|K_k(x_B,x-y)\psi_0(2^{-\ell}(x-y)) -  K_k(x_B,x-y)\psi_0(2^{-\ell}(x_B-y))\right|^{p'} dy \right)^\frac{1}{p'}
\end{equation*}
(where $K_k$ has the same meaning as in \cite{Michalowski_Rule_Staubach_2012}) needs to be controlled, in addition to $I_{3,1}$ and $I_{3,2}$. However, $I_{3,3}$ can be estimated as in (3.6) in \cite{Michalowski_Rule_Staubach_2012}, as we only need that derivatives of $\psi_0$ are $O(1) \leq O(2^k)$.
\end{proof}

We now prove Theorem~\ref{sharpestimatethm}, which expands on the results of Chanillo \& Torchinsky~\cite{ChanilloTorchinsky1986}, Alvarez \& Hounie~\cite{AlvarezHounie90}, Chen \& Wang~\cite{ChenWang2023} and Wang~\cite{wang2022sharpfunctionweightedlp}.

\begin{proof}[Proof of Theorem~\ref{sharpestimatethm}]
The method of proof is essentially contained in the proof of Theorem~3.2 in \cite{AlvarezHounie90}. To prove \eqref{sharppwbound1}, it suffices to prove the estimate
\begin{equation*}
    \frac{1}{|B|} \int_{B} | T_a(f)(x) - ( T_a(f) )_B| dx \lesssim M_p(f)(x_0)
\end{equation*}
for all balls $B$ containing $x_0$. Let $\tau$ denote the radius of $B$. We begin by letting $B'$ be the ball concentric to $B$ with radius $2\tau^\rho$. We decompose $f$ as
\begin{equation*}
f= f\chi_{B'}+ f(1-\chi_{B'})=: f_1 +f_2
\end{equation*}
and so
\begin{equation}\label{eq:iandii}
\begin{aligned}
&\frac{1}{|B|} \int_{B}|T_a(f)(x) - ( T_a(f))_{B}| dx \\
&\leq \frac{2}{|B|} \int_{B} |T_a(f_1)(x)| dx + \frac{1}{|B|} \int_{B} |T_a(f_2)(x) - ( T_a(f_2) )_{B}| dx\\
&:= \textbf{I} +\textbf{II}.
\end{aligned}
\end{equation}

To estimate $\textbf{I}$, we can make use of the fact that $T_a$ is bounded from $L^{p}$ to $L^{p/\rho}$ when $1 < p \leq 2 \leq p/\rho$ (see Theorem~3.5~(a) in \cite{AlvarezHounie90} or even Lemma~3.1 there if $p=2$). Thus,
\begin{equation*}
    \textbf{I} \leq \left(\frac{2}{|B|} \int_{B} |T_a(f_1)(x)|^{p/\rho} dx\right)^{\frac{\rho}{p}}
    \lesssim \frac{2}{|B|^{\rho/p}}\left( \int |f_1(x)|^{p} dx\right)^{\frac{1}{p}}
    \lesssim \left\langle f \right\rangle_{p,B'} \leq M_p(f)(x_0).
\end{equation*}

To estimate $\textbf{II}$ we begin with the case $\tau \geq 1$. We can use the kernel estimate
\begin{equation*}
    |K(x,z)| \lesssim |z|^{-N} \quad \mbox{for $N>0$, $|z| \geq 1$ and $\rho > 1$}
\end{equation*}
from, for example, (3.2) in \cite{Michalowski_Rule_Staubach_2012}. This gives us the estimate
\begin{equation}\label{largecube}
\begin{aligned}
    |T_a(f_2)(x)| &\leq \int |K(x,x-y)||f_2(y)| dy \\
    &\lesssim \sum_{j=1}^\infty \int_{A_j} |x-y|^{-(n+1)}|f_2(y)| dy \\
    &\lesssim \sum_{j=1}^\infty \frac{1}{2^{j(n+1)}r^{(n+1)\rho}}\int_{A_j}|f_2(y)| dy \\
    &\lesssim \sum_{j=1}^\infty \frac{1}{2^jr^{\rho}}Mf_2(x_0)
    \lesssim Mf_2(x_0)
\end{aligned}
\end{equation}
for $x \in B$, where $A_j = \{y \,|\, 2^j\tau^\rho \leq |y-x_B| \leq 2^{j+1}\tau^\rho\}$ and $x_B$ is the common centre of $B$ and $B'$. 

Switching now to the case $\tau < 1$ we can apply Lemma~\ref{kernelestimatelemma} with $c_1= c_2 =1$ and $\theta =\rho$. Denoting the centre of $B$ again by $x_B$ and applying \eqref{kernellem} we obtain
\begin{equation}
\begin{aligned}\label{smallcube}
&|T_a(f_2)(x) - T_a(f_2)(x_B)|\\
&\leq \left| \int (K(x, x -y)- K(x_B, x_B -y)) f_2 (y) dy\right| \\
&\leq \sum_{j=1}^{\infty} \left(\int_{A_j}| K(x, x -y)- K(x_B, x_B -y)|^{p'}dy\right)^{\frac{1}{p'}}\left(\int_{A_j} | f_2 (y)|^{p} dy\right)^{\frac{1}{p}}\\
&\lesssim \sum_{j=1}^{\infty}  2^{-jh} \tau^{-m-\frac{n}{p}} (2^{j} \tau^\rho)^{\frac{n}{p}} M_{p}f_2(x_0)
\lesssim \sum_{j=1}^{\infty}  2^{j(\frac{n}{p}-h)}\tau^{\frac{n}{p}(\rho - 1) - m} M_{p}f_2(x_0),
\end{aligned}
\end{equation}
for $x\in B$ and an appropriate choice of $h$. The conditions $h$ needs to satisfy are $h>\frac{n}{p}$, which ensures the series $\sum_{j\geq1}2^{j(\frac{n}{p}-l)}$ is finite, and those of Lemma~\ref{kernelestimatelemma}. Since
\begin{equation*}
    m + \frac{n}{p} = \frac{n\rho}{p} - n\lambda \leq \frac{n\rho}{p},
\end{equation*}
the conditions of the lemma and $h > \frac{n}{p}$ can be combined to read
\begin{equation*}
    \frac{n}{p} < h < \frac{n}{p} + \frac{1-n\lambda}{\rho}.
\end{equation*}
Clearly such an $h$ exists under the condition that $\lambda < 1/n$.

Now we can estimate $\textbf{II}$ using \eqref{largecube} or \eqref{smallcube}. We write
\begin{align*}
     \textbf{II} &= \frac{1}{|B|} \int_{B} |T_a(f_2)(x) - ( T_a(f_2) )_{B}| dx \\
     &= \frac{1}{|B|} \int_{B} \left|\frac{1}{|B|} \int_{B}\left(T_a(f_2)(x) - T_a(f_2)(z)\right)dz\right| dx \\
     &\leq \frac{1}{|B|^2} \int_{B} \int_{B} \left|T_a(f_2)(x) - T_a(f_2)(z)\right| dz dx,
\end{align*}
so for $\tau \geq 1$ we use \eqref{largecube} to estimate
\begin{equation*}
    \textbf{II} \leq \sup_{x\in B} |T_a(f_2)(x)| + \sup_{z\in B} |T_a(f_2)(z)| \lesssim Mf_2(x_0) \leq M_pf(x_0)
\end{equation*}
and for $\tau < 1$ we use \eqref{smallcube} to estimate
\begin{equation*}
    \textbf{II} \leq \sup_{x\in B} |T_a(f_2)(x) - T_a(f_2)(x_B)| +  \sup_{z\in B} |T_a(f_2)(x_B) - T_a(f_2)(z)|
     \lesssim M_pf_2(x_0) \leq M_pf(x_0).
\end{equation*}
Since \eqref{largecube} and \eqref{smallcube} apply to any $x\in B$, they also apply for $x = z$.
\end{proof}

We are now ready to piece together the estimates needed to apply Proposition~\ref{intro:end-pointsparsebdd2} and prove Theorem~\ref{thm:endpoint-pseudos}. Recall that $\psi_0$ is defined on page~\pageref{page:psi_0}.

\begin{lem}\label{off-diagonal_extremes}
Fix $\eta \in (0,1)$, $0\leq \delta < 1$ and suppose $a\in S^{0}_{1,\delta}$. Define
\begin{equation}\label{amplitude}
    \tilde{a}_{\ell}(x,y,\xi) = a(x,\xi)\psi_0(2^{-\ell}(x-y)).
\end{equation}
Then, for each bounded and compactly supported $f$, there exists a $(3^{-n}\eta)$-sparse collection of cubes $\mathcal{S}$ which forms a graded partially ordered set with partial order $\preceq$ and rank function $\varrho$ satisfying \ref{assump:inductionaxiom}, \ref{assump:flag} and \ref{assump:support} from Section~\ref{sec:2} with $T = T_{a_\ell}$. Moreover if $Q',Q \in \mathcal{S}$ and $Q' \preceqdot Q$ we have that
\begin{equation}
    \sup_{x\in \frac{1}{3}Q'}\left|T_{\tilde{a}_\ell}(f\chi_{Q\setminus Q'})(x) - ( T_{\tilde{a}_\ell}(f\chi_{Q\setminus Q'}) )_{\frac{1}{3}Q'}\right| \lesssim \langle f \rangle_{1,Q},\label{bmointerpolate}
\end{equation}
uniformly in $\ell$. Furthermore, $\mathcal{S}$ can be chosen so that, for a given bounded and compactly supported function $g$ and for given $1 < r \leq s \leq \infty$, we have
\begin{align}
        \left\langle \tilde{g}_Q \right\rangle_{r',Q} &\lesssim \left\langle g \right\rangle_{s',Q}
        \quad \mbox{and} \label{greverseholder} \\
        \left\langle g \right\rangle_{s',Q'} &\lesssim \left\langle g \right\rangle_{s',Q}, \label{gsmalltolarge}
\end{align}
where $\tilde{g}_Q$ is defined in \eqref{def:gtildeq} and again $Q' \preceqdot Q$.
\end{lem}

\begin{proof}
We construct our sparse set in a similar manner to that used in the proof of Proposition~\ref{intro:sparseform} in Section~\ref{page:constructionofcubes}. It will be convenient to first construct
\begin{equation}
    \mathcal{R} := \bigcup_{p=1}^\infty \mathcal{R}_p \quad \mbox{and} \quad \mathcal{S} := \bigcup_{p=1}^\infty \mathcal{S}_p, \label{def:secondsparseset}
\end{equation}
where each $\mathcal{R}_p$ is defined recursively in $p$, and then define $\mathcal{S}_p$ to be the set of concentric triples $3Q$ of cubes $Q \in \mathcal{R}_p$. The rank $\varrho(Q)$ of a cube $Q$ will be defined to be $p$ when $Q \in \mathcal{S}_p$. To do this, we first select a collection of cubes $\mathcal{R}_0 := \{Q^{m,0}_k\}_m$ (with the notation defined in \eqref{def:dyadic_cube}) whose union contains the support of $f$ and $g$ and with sufficiently large side length (that is, sufficiently small $k$) so that the supports of $T_{\tilde{a}_{\ell}}(f\chi_{Q^{m,0}_k})$ are contained in the concentric triple $3Q^{m,0}_k$. This is possible due to the support properties of the kernel of $ T_{\tilde{a}_{\ell}}$. Thus, if we succeed in completing this recursive process to define $\mathcal{S}$, it will satisfy \ref{assump:support}.

To carry out this process, for a given $Q \in \mathcal{R}_{p-1}$, we will define $\mathcal{R}_p$ by defining a collection $\mathcal{R}(Q)$ of sub-cubes of $Q$ and then define
\begin{equation*}
    \mathcal{R}_{p} = \bigcup_{Q \in \mathcal{R}^{\omega}_{p-1}} \mathcal{R}(Q).
\end{equation*}
Consider the sets 
\begin{equation}\label{def:stoppingtimecubes}
\begin{aligned}
    F^c_1(Q) &= \left\{ x \in \R^n : M^c_r(f\chi_{3Q})(x) > \left(\frac{2\cdot 3^n}{1-\eta}\right)^{1/r}\|M^c_{r}\|_{L^r\to L^{r,\infty}}\langle f \rangle_{r,3Q} \right\} \quad \mbox{and} \\
    F^c_2(Q) &= \left\{ x \in \R^n : M^c_{s'}(g\chi_Q)(x) > \left(\frac{2}{1-\eta}\right)^{1/s'}\|M^c_{s'}\|_{L^{s'}\to L^{s',\infty}}\langle g \rangle_{s',Q} \right\}.
\end{aligned}
\end{equation}
For each $j=1,2$, let $\tilde{\mathcal{R}}^j(Q)$ be a collection of dyadic cubes which form a Whitney decomposition of $F^c_j(Q)$ (see Appendix J in \cite{grafakos2014classical}). The weak-type $(r,r)$ boundedness of $M^c_r$ means that
\begin{align*}
    |F^c_1(Q)| &\leq \left(\frac{1-\eta}{2\cdot 3^n\|M^c_r\|_{L^r\to L^{r,\infty}}^2\langle f \rangle_{r,3Q}^r}\right)\|M^c_r\|_{L^r\to L^{r,\infty}}^r \int |f\chi_{3Q}|^r \\
    &\leq \left(\frac{1-\eta}{2\cdot 3^n\langle f \rangle_{r,3Q}^r}\right) \langle f \rangle_{r,3Q}^r|3Q|
    \leq \left(\frac{1-\eta}{2}\right) |Q|,
\end{align*}
and a similar reasoning shows that
\begin{equation*}
    |F^c_2(Q)| \leq \left(\frac{1-\eta}{2}\right) |Q|.
\end{equation*}
Therefore, $|F^c_1(Q) \cup F^c_2(Q)| \leq (1-\eta)|Q|$ and, in particular, any cube in the Whitney decompositions of $F^c_j(Q)$ is either contained in $Q$ or disjoint from it. Consequently, we can set
\begin{equation*}
    \tilde{\mathcal{R}}(Q) := \left\{ Q' \in \tilde{\mathcal{R}}^1(Q)\cup\tilde{\mathcal{R}}^2(Q)\colon Q' \subseteq Q\right\}
\end{equation*}
to be the Whitney cubes which are also contained in $Q$, and define
\begin{equation*}
    \mathcal{R}(Q) = \{Q' \in \tilde{\mathcal{R}}(Q) \colon \mbox{$Q' \not\subset Q''$ for any other $Q'' \in \tilde{\mathcal{R}}(Q)$}\}
\end{equation*}
to be the collection of maximal cubes of $\tilde{\mathcal{R}}(Q)$ with respect to inclusion. We define
\begin{equation*}
    F^c(Q) = F^c_1(Q) \cup F^c_2(Q),
\end{equation*}
and $F(Q) = Q\setminus F^c(Q)$ for each dyadic cube $Q \in \mathcal{R}_{p-1}$. Thus, we have decomposed
\begin{equation*}\label{coveringrelation}
    Q = F(Q) \cup \left(\bigcup_{Q' \in \mathcal{R}(Q)} Q'\right)
\end{equation*}
as a disjoint union of $F(Q)$ and the disjoint dyadic cubes $Q' \in \mathcal{R}(Q)$, and $|F(Q)| \geq (1-(1-\eta)) |Q| \geq \eta|Q|$. Thus, we have shown that $\mathcal{R}$ is $\eta$-sparse, and consequently $\mathcal{S}$ defined in the paragraph containing \eqref{def:secondsparseset} is $(3^{-n}\eta)$-sparse. Observe that $\mathcal{R}(Q)$ is the set of all cubes in $\mathcal{R}$ which are covered by $Q$ with respect to the partial ordering of inclusion, so $\mathcal{S}$ inherits the partial ordering \eqref{finclusion} from $\mathcal{R}$. We also see that $F(Q) = E(3Q)$ with $E(Q)$ defined in \eqref{def:f}. We conclude that $\mathcal{S}$ satisfies \ref{assump:inductionaxiom}. It is also clear from the iterative construction that the rank function satisfies \ref{assump:flag}.

Observe that the Whitney decomposition is such that $4\sqrt{n}Q'$ intersects the compliment of $E^c_j(Q)$ for each $j=1,2$ and $Q' \in \mathcal{R}(Q)$ and consequently there exist an $x$ in this intersection such that the reverse of the inequalities in \eqref{def:stoppingtimecubes} holds. We will use this fact to prove \eqref{bmointerpolate}--\eqref{gsmalltolarge}.

To prove \eqref{bmointerpolate} we first recall the proof of Theorem~\ref{sharpestimatethm}. Looking at the proof of $\textbf{II}$ in \eqref{eq:iandii}, we see that we actually estimate the integrand of $\textbf{II}$, not just the integral. So, using the same notation, we have the estimate
\begin{equation*}
   \sup_{x\in B} \left|T_{a}(f_2)(x) - ( T_{a}(f_2) )_{B}\right| \lesssim M(f_2)(x_0)
\end{equation*}
where the support of $f_2$ does not intersect $B'$ and $x_0 \in B$. Noting that \eqref{kernellem2} lets us run the same estimates of $\textbf{II}$ in the proof of Theorem~\ref{sharpestimatethm} with $T_a$ replaced by $T_{a_\ell}$ and, when $\rho=1$, $B'$ is the concentric double of $B$, we can conclude that
\begin{equation*}
    \sup_{x\in Q'}\left|T_{\tilde{a}_{\ell}}(f\chi_{3Q\setminus 3Q'})(x) - ( T_{\tilde{a}_{\ell}}(f\chi_{3Q\setminus 3Q'}) )_{Q'}\right| \lesssim M(f\chi_{3Q\setminus 3Q'})(x_0)
\end{equation*}
when $a\in S^0_{1,\delta}$ and $x_0 \in Q'$. However, observe that if $x_0 \in Q'$ and $B$ is a ball containing $x_0$ and of radius less than half the side length of $Q'$, then the average of $f\chi_{3Q\setminus 3Q'}$ over $B$ is zero. For balls $B$ containing $x_0 \in Q'$ with radius at least half the side length of $Q'$, then for any $\tilde{x} \in 4\sqrt{n}Q'$, the ball $\tilde{B}$ centred at $\tilde{x}$ with radius $8\sqrt{n}$ times the radius of $B$ will contain $B$. Consequently
\begin{equation*}
    M_r(f\chi_{3Q\setminus 3Q'})(x_0) \lesssim \inf_{\tilde{x} \in 4\sqrt{n}Q'}M_r(f\chi_{3Q\setminus 3Q'})(\tilde{x}) \leq \inf_{\tilde{x} \in 4\sqrt{n}Q'}M^c_r(f\chi_{3Q})(\tilde{x}) \lesssim \langle f \rangle_{r,3Q}
\end{equation*}
by the reverse of the first inequality in \eqref{def:stoppingtimecubes}. Rescaling the cubes by a factor $\frac{1}{3}$ proves \eqref{bmointerpolate}.

Since $F(Q) \subseteq Q \setminus F^c_2(Q)$, we know that
\begin{equation*}
    g(x_0) \leq M^c_{s'}(g\chi_Q)(x_0) \leq \left(\frac{4}{1-\eta}\right)^{1/s'}\|M^c_{s'}\|_{L^{s'}\to L^{s',\infty}}\langle g \rangle_{s',Q}
\end{equation*}
for all $x_0 \in F(Q)$. Moreover, for each $Q' \in \mathcal{R}(Q)$, we know that $F(Q) \cap 3Q' \neq \emptyset$, so there exists an $x_0 \in F(Q) \cap 3Q'$ such that
\begin{equation*}
    |g_{Q'}| \leq \langle g\chi_Q \rangle_{s',3Q'} \leq M^c_{s'}(g\chi_Q)(x_0) \leq \left(\frac{4}{1-\eta}\right)^{1/s'}\|M^c_{s'}\|_{L^{s'}\to L^{s',\infty}}\langle g \rangle_{s',Q}.
\end{equation*}
These last two estimates give us that
\begin{equation*}
    \langle \tilde{g}_Q \rangle_{\infty,Q} \lesssim \langle g \rangle_{s',Q}
\end{equation*}
and we obtain \eqref{greverseholder} from this.

We can prove \eqref{gsmalltolarge} in a similar manner. By construction, we know that when $Q'\in\mathcal{R}(Q)$, there exists an $x_0 \in 4\sqrt{n}Q'$ such that
\begin{equation*}
    M^c_{s'}(g\chi_Q)(x_0) \leq \left(\frac{4}{1-\eta}\right)^{1/s'}\|M^c_{s'}\|_{L^{s'}\to L^{s',\infty}}\langle g \rangle_{s',Q}
\end{equation*}
and, since
\begin{equation*}
    \langle g \rangle_{s',Q'} \lesssim \langle g\chi_Q \rangle_{s',4\sqrt{n}Q'} \leq M^c_{s'}(g\chi_Q)(x_0),
\end{equation*}
\eqref{gsmalltolarge} follows.
\end{proof}

Finally, we state and prove an application of Proposition~\ref{intro:end-pointsparsebdd2}.
\begin{thm}\label{thm:endpoint-pseudos}
    Suppose $a\in S^0_{1,\delta}$ with $\delta \in [0,1)$ and $r > 1$. Then there exists a constant $C$, which only depends on $\delta$, $r$ and a finite number of the implicity constants in Defintion~\ref{def:symbol} such that for each pair of bounded functions $f$ and $g$ with compact support, there exists a sparse collection $\mathcal{S}$ such that
\begin{equation*}
    \left|\langle T_a(f),g\rangle\right| \leq C \sum_{Q\in\mathcal{S}} \langle f \rangle_{r,Q}\langle g \rangle_{1,Q}|Q|.
\end{equation*}
\end{thm}
Observe that this theorem also yields a weak-type result of the kind discussed in Remark~\ref{commentc} when $m=0$ and $\rho=1$. Indeed, \eqref{ineq:weak-type} can be seen to hold when $a \in S^0_{1,\delta}$ with $\delta \in [0,1)$ by applying the fact that these symbols are closed under adjoints and Theorem~\ref{thm:endpoint-pseudos} to Theorem~E in \cite{Conde-AlonsoCuliucPlinioOu}. However, in the course of proving Theorem~\ref{thm:endpoint-pseudos}, we will use the kernel estimates of Lemma~\ref{kernelestimatelemma}, which are very similar to the estimates used in Theorem~3.2 of \cite{AlvarezHounie90} to give a direct proof of weak-type boundedness.

\begin{proof}[Proof of Theorem~\ref{thm:endpoint-pseudos}]
We will prove Theorem~\ref{thm:endpoint-pseudos} by applying Proposition~\ref{intro:end-pointsparsebdd2} to the operator $T = T_{\tilde{a}_{\ell}}$ with $\tilde{a}_\ell$ as in \eqref{amplitude}. In order to justify the assumptions of Proposition~\ref{intro:end-pointsparsebdd2} we invoke Lemma~\ref{off-diagonal_extremes}. This proves the existence of a sparse set $\mathcal{S}$ which satisfies assumptions \ref{assump:inductionaxiom}, \ref{assump:flag} and \ref{assump:support}. Moreover, \eqref{bmointerpolate}, \eqref{greverseholder} and \eqref{gsmalltolarge} prove \ref{assump:lrls} with $s=\infty$ and $r\geq1$, \eqref{end-point2} and \eqref{end-point3}, respectively. The operator arising from the amplitude in \eqref{amplitude} can be written as an operator arising from a symbol in $S^0_{1,\delta}$ (see the proposition on page~258 of \cite{S}). Thus, we can apply Theorem~3.4 in \cite{AlvarezHounie90} to conclude the operator is bounded from $L^r$ to $L^r$ for any $r>1$ and we conclude \eqref{end-point1} holds, with $A_2$ independent of $\ell$.

Thus, all the hypotheses of Proposition~\ref{intro:end-pointsparsebdd2} hold, and we conclude that
\begin{equation}\label{ineq:sparselocalpart}
    \left|\langle T_{\tilde{a}_{\ell}}(f),g\rangle\right| \lesssim \sum_{Q\in\mathcal{S}} |Q|\langle f \rangle_{r,Q}\langle g \rangle_{1,Q},
\end{equation}
with an implicit constant uniform in $\ell$.

From here we have two alternative ways to complete the proof. We can either apply Lemma~4.7 in \cite{LaceyMena} to replace $\mathcal{S}$ in \eqref{ineq:sparselocalpart} with a maximal sparse collection for $f$ and $g$. Then taking the supremum in $\ell$ proves Theorem~\ref{thm:endpoint-pseudos}. Or the kernel estimates in Chapter VI, \S\S~2.2-2.3 of \cite{S} can be used to prove a sparse bound for $T_{a-\tilde{a}_\ell}$. Thus, the operator $\langle T_a(f), g \rangle$ is bounded by the sum of two sparse bounds, which by Theorem~1.3 of \cite{Hanninen} is itself a sparse bound and so the proof of Theorem~\ref{thm:endpoint-pseudos} is again complete.
\end{proof}

%%%%%%%%%%%%%%%%%%%%%%
%%%%%%%%%%%%%%%%%%%%%%
%%%%%%%%%%%%%%%%%%%%%%
%%%%%%%%%%%%%%%%%%%%%%

\subsubsection*{Acknowledgements} We would like to thank the anonymous referee for many helpful suggestions which improved this paper, and, in particular, for pointing out Remark~\ref{commentc} to us. We would also like to thank the Swedish International Development Cooperation Agency for funding this work.

\printbibliography

@book {S,
    AUTHOR = {Stein, Elias M.},
     TITLE = {Harmonic analysis: real-variable methods, orthogonality, and
              oscillatory integrals},
    SERIES = {Princeton Mathematical Series},
    VOLUME = {43},
      NOTE = {With the assistance of Timothy S. Murphy,
              Monographs in Harmonic Analysis, III},
 PUBLISHER = {Princeton University Press, Princeton, NJ},
      YEAR = {1993},
     PAGES = {xiv+695},
      ISBN = {0-691-03216-5},
   MRCLASS = {42-02 (35Sxx 43-02 47G30)},
  MRNUMBER = {1232192},
MRREVIEWER = {Michael\ Cowling},
}

@article {MRS1,
    AUTHOR = {Michalowski, Nicholas and Rule, David and Staubach,
              Wolfgang},
     TITLE = {Weighted norm inequalities for pseudo-pseudodifferential
              operators defined by amplitudes},
   JOURNAL = {J. Funct. Anal.},
  FJOURNAL = {Journal of Functional Analysis},
    VOLUME = {258},
      YEAR = {2010},
    NUMBER = {12},
     PAGES = {4183--4209},
      ISSN = {0022-1236},
   MRCLASS = {35S05 (42B25 47G30)},
  MRNUMBER = {2609542},
       DOI = {10.1016/j.jfa.2010.03.013},
       URL = {https://doi.org/10.1016/j.jfa.2010.03.013},
}

@book {grafakos2014classical,
    AUTHOR = {Grafakos, Loukas},
     TITLE = {Classical {F}ourier analysis},
    SERIES = {Graduate Texts in Mathematics},
    VOLUME = {249},
   EDITION = {Third},
 PUBLISHER = {Springer, New York},
      YEAR = {2014},
     PAGES = {xviii+638},
      ISBN = {978-1-4939-1193-6},
   MRCLASS = {42-01 (42Bxx)},
  MRNUMBER = {3243734},
MRREVIEWER = {Atanas\ G.\ Stefanov},
       DOI = {10.1007/978-1-4939-1194-3},
       URL = {https://doi.org/10.1007/978-1-4939-1194-3},
}

@article {BC,
    AUTHOR = {Beltran, David and Cladek, Laura},
     TITLE = {Sparse bounds for pseudodifferential operators},
   JOURNAL = {J. Anal. Math.},
  FJOURNAL = {Journal d'Analyse Math\'ematique},
    VOLUME = {140},
      YEAR = {2020},
    NUMBER = {1},
     PAGES = {89--116},
      ISSN = {0021-7670,1565-8538},
   MRCLASS = {42B25 (35R11 35S05)},
  MRNUMBER = {4094458},
MRREVIEWER = {Yaoming\ Niu},
       DOI = {10.1007/s11854-020-0083-x},
       URL = {https://doi.org/10.1007/s11854-020-0083-x},
}

@article{Michalowski_Rule_Staubach_2012, title={Weighted Lp Boundedness of Pseudodifferential Operators and Applications}, volume={55}, DOI={10.4153/CMB-2011-122-7}, number={3}, journal={Canadian Mathematical Bulletin}, author={Michalowski, Nicholas and Rule, David J. and Staubach, Wolfgang}, year={2012}, pages={555–570}}

@article {AlvarezHounie90,
    AUTHOR = {\'Alvarez, Josefina and Hounie, Jorge},
     TITLE = {Estimates for the kernel and continuity properties of
              pseudo-differential operators},
   JOURNAL = {Ark. Mat.},
  FJOURNAL = {Arkiv f\"or Matematik},
    VOLUME = {28},
      YEAR = {1990},
    NUMBER = {1},
     PAGES = {1--22},
      ISSN = {0004-2080,1871-2487},
   MRCLASS = {35S05 (42B30 47G30)},
  MRNUMBER = {1049640},
MRREVIEWER = {Luigi\ Rodino},
       DOI = {10.1007/BF02387364},
       URL = {https://doi.org/10.1007/BF02387364},
}

@article {ChenWang2023,
    AUTHOR = {Chen, Wenyi and Wang, Guangqing},
     TITLE = {A pointwise estimate for pseudo-differential operators},
   JOURNAL = {Bull. Math. Sci.},
  FJOURNAL = {Bulletin of Mathematical Sciences},
    VOLUME = {13},
      YEAR = {2023},
    NUMBER = {2},
     PAGES = {Paper No. 2250001, 13},
      ISSN = {1664-3607,1664-3615},
   MRCLASS = {42B20 (42B37)},
  MRNUMBER = {4641743},
MRREVIEWER = {Justin\ G.\ Trulen},
       DOI = {10.1142/S1664360722500011},
       URL = {https://doi.org/10.1142/S1664360722500011},
}

@misc{wang2022sharpfunctionweightedlp,
      title={Sharp function and weighted $L^{p}$ estimates for pseudo-differential operators with symbols in general H\"{o}rmander classes}, 
      author={Guangqing Wang},
      year={2022},
      eprint={2206.09825},
      archivePrefix={arXiv},
      primaryClass={math.AP},
      url={https://arxiv.org/abs/2206.09825}, 
}

@article {ChanilloTorchinsky1986,
    AUTHOR = {Chanillo, Sagun and Torchinsky, Alberto},
     TITLE = {Sharp function and weighted {$L^p$} estimates for a class of
              pseudodifferential operators},
   JOURNAL = {Ark. Mat.},
  FJOURNAL = {Arkiv f\"or Matematik},
    VOLUME = {24},
      YEAR = {1986},
    NUMBER = {1},
     PAGES = {1--25},
      ISSN = {0004-2080,1871-2487},
   MRCLASS = {47G05 (35S05)},
  MRNUMBER = {852824},
MRREVIEWER = {Yves\ Meyer},
       DOI = {10.1007/BF02384387},
       URL = {https://doi.org/10.1007/BF02384387},
}

@article {BernicotFreyPetermichl,
    AUTHOR = {Bernicot, Fr\'ed\'eric and Frey, Dorothee and Petermichl,
              Stefanie},
     TITLE = {Sharp weighted norm estimates beyond {C}alder\'on-{Z}ygmund
              theory},
   JOURNAL = {Anal. PDE},
  FJOURNAL = {Analysis \& PDE},
    VOLUME = {9},
      YEAR = {2016},
    NUMBER = {5},
     PAGES = {1079--1113},
      ISSN = {2157-5045,1948-206X},
   MRCLASS = {42B20 (58J35)},
  MRNUMBER = {3531367},
MRREVIEWER = {Elena\ Cordero},
       DOI = {10.2140/apde.2016.9.1079},
       URL = {https://doi.org/10.2140/apde.2016.9.1079},
}

@article {Conde-AlonsoCuliucPlinioOu,
    AUTHOR = {Conde-Alonso, Jos\'e{} M. and Culiuc, Amalia and Di Plinio,
              Francesco and Ou, Yumeng},
     TITLE = {A sparse domination principle for rough singular integrals},
   JOURNAL = {Anal. PDE},
  FJOURNAL = {Analysis \& PDE},
    VOLUME = {10},
      YEAR = {2017},
    NUMBER = {5},
     PAGES = {1255--1284},
      ISSN = {2157-5045,1948-206X},
   MRCLASS = {42B20 (42B25)},
  MRNUMBER = {3668591},
MRREVIEWER = {Liguang\ Liu},
       DOI = {10.2140/apde.2017.10.1255},
       URL = {https://doi.org/10.2140/apde.2017.10.1255},
}

@article {LiPerezRivera-RiosRoncal,
    AUTHOR = {Li, Kangwei and P\'erez, Carlos and Rivera-R\'ios, Israel P.
              and Roncal, Luz},
     TITLE = {Weighted norm inequalities for rough singular integral
              operators},
   JOURNAL = {J. Geom. Anal.},
  FJOURNAL = {Journal of Geometric Analysis},
    VOLUME = {29},
      YEAR = {2019},
    NUMBER = {3},
     PAGES = {2526--2564},
      ISSN = {1050-6926,1559-002X},
   MRCLASS = {42B20 (42B15 42B25)},
  MRNUMBER = {3969435},
MRREVIEWER = {Luc\ Del\'eaval},
       DOI = {10.1007/s12220-018-0085-4},
       URL = {https://doi.org/10.1007/s12220-018-0085-4},
}

@online{LMSvideos,
        title = {Lecture 2: Sparse domination in harmonic analysis},
        date = {2021},
        organization = {London Mathematical Society},
        author = {Conde-Alonso, Jos\'e{} M.},
        url = {https://www.youtube.com/watch?v=EuH_uq2qrEY},
    }

@article {Hanninen,
    AUTHOR = {H\"anninen, Timo S.},
     TITLE = {Equivalence of sparse and {C}arleson coefficients for general
              sets},
   JOURNAL = {Ark. Mat.},
  FJOURNAL = {Arkiv f\"or Matematik},
    VOLUME = {56},
      YEAR = {2018},
    NUMBER = {2},
     PAGES = {333--339},
      ISSN = {0004-2080,1871-2487},
   MRCLASS = {42B15 (46E35)},
  MRNUMBER = {3893778},
MRREVIEWER = {Oscar\ Blasco},
       DOI = {10.4310/ARKIV.2018.v56.n2.a8},
       URL = {https://doi.org/10.4310/ARKIV.2018.v56.n2.a8},
}

@article{Lerner2016,
  author    = {Andrei K. Lerner},
  title     = {On pointwise estimates involving sparse operators},
  journal   = {New York J. Math.},
  volume    = {22},
  year      = {2016},
  pages     = {341--349},
  eprint    = {1603.01330},
  archivePrefix = {arXiv},
  primaryClass = {math.CA}
}

@article {Hounie,
    AUTHOR = {Hounie, Jorge},
     TITLE = {On the {$L^2$} continuity of pseudodifferential operators},
   JOURNAL = {Comm. Partial Differential Equations},
  FJOURNAL = {Communications in Partial Differential Equations},
    VOLUME = {11},
      YEAR = {1986},
    NUMBER = {7},
     PAGES = {765--778},
      ISSN = {0360-5302,1532-4133},
   MRCLASS = {35S05 (47G05)},
  MRNUMBER = {837930},
MRREVIEWER = {Yves\ Meyer},
       DOI = {10.1080/03605308608820444},
       URL = {https://doi.org/10.1080/03605308608820444},
}

@article {LaceyMena,
    AUTHOR = {Lacey, Michael T. and Mena Arias, Dar\'io},
     TITLE = {The sparse {T}1 theorem},
   JOURNAL = {Houston J. Math.},
  FJOURNAL = {Houston Journal of Mathematics},
    VOLUME = {43},
      YEAR = {2017},
    NUMBER = {1},
     PAGES = {111--127},
      ISSN = {0362-1588},
   MRCLASS = {42B20 (42B25)},
  MRNUMBER = {3647935},
MRREVIEWER = {Joan\ Orobitg},
}

@phdthesis{BeltranPortales,
  title        = {Geometric control of oscillatory integrals},
  author       = {David Beltran Portalés},
  year         = 2017,
  note         = {Available at \url{https://etheses.bham.ac.uk/id/eprint/7566/}},
  school       = {University of Birmingham},
  type         = {PhD thesis}
}

\Addresses

\end{document}